\documentclass[10pt]{article}
\usepackage{graphicx,epsfig}
\usepackage{amsmath}
\usepackage[utf8]{inputenc}
\usepackage{newunicodechar}
\newunicodechar{　}{~} 

 \usepackage{color}      

\usepackage{url}

\usepackage{tabu}

\usepackage{array}
\usepackage{makecell}

\usepackage{multirow}

\usepackage{subcaption}
\usepackage[driverfallback=dvipdfm,
bookmarks=true,    
bookmarksopen=true,
bookmarksopenlevel=\maxdimen,
bookmarksdepth=10
]{hyperref}
\hypersetup{
pdfcenterwindow=true, 
pdfpagelabels=true,
pagebackref=true,            
hypertexnames=true,
plainpages=false,
unicode=false,     
pdftoolbar=true,                
    pdfmenubar=true,          
    pdffitwindow=false,         
    pdfstartview={FitH},        
    pdftitle={On Robust Alpha-Damping Viscous Scheme},
    pdfauthor={Hiroaki Nishikawa},    
    pdfsubject={On Robust Alpha-Damping Viscous Scheme},       
    pdfcreator={Hiroaki Nishikawa},   
    pdfproducer={Hiroaki Nishikawa}, 
    baseurl={https://www.researchgate.net/profile/Hiroaki-Nishikawa-2},
    pdfkeywords={}, 
    pdfnewwindow=true,     
    colorlinks=true,       
    linkcolor=black,       
    citecolor=black,       
    filecolor=black,       
    urlcolor=black,        
  breaklinks=true,
  hyperfigures=true,
  backref=true,
breaklinks=false, 
pdfpagelabels,      
pagebackref,        
hypertexnames=true, 
plainpages=false,   
naturalnames        
}
\usepackage[all]{hypcap}

\numberwithin{equation}{section}
\numberwithin{subsubsection}{subsection}
\numberwithin{subsection}{section}

\usepackage[thinlines]{easytable}

\usepackage{cleveref}
\usepackage[thinlines]{easytable}

\usepackage{blindtext,graphicx}
\usepackage[absolute]{textpos}
\numberwithin{equation}{section}
\numberwithin{subsubsection}{subsection}
\numberwithin{subsection}{section}

\definecolor{black}{rgb}{0.43, 0.21, 0.1} 
\usepackage{enumitem,xcolor}

\begin{document}


\title{\bf On Robust Alpha-Damping Viscous Scheme}
 
 \author{
{Hiroaki Nishikawa} \\
  {\itshape {National Institute of Aerospace}, Hampton, VA 23666, USA}
}

\date{} 

\maketitle 


\begin{abstract}
In this paper, we investigate the convergence of an implicit defect-correction solver for a viscous discretization based on the alpha-damping scheme for unstructured grids. We show that significantly more robust iterative convergence is achieved by evaluating the damping term at the midpoint between two adjacent cell centers (edge midpoint) rather than at the face centroid. A one-dimensional Fourier analysis reveals that the implicit solver tends to be stable when the damping term in the residual is smaller than that used to construct the Jacobian. 
This observation suggests that the solver can be stabilized by effectively reducing the magnitude of the damping term in the residual — an effect achieved by the edge-midpoint evaluation. Robust convergence is demonstrated numerically for two-dimensional viscous-flow problems on highly irregular mixed-element and triangular grids. 
\end{abstract}

\section{Introduction}
\label{intro} 

\indent

Second-order finite-volume methods are widely used in practical computational fluid dynamics (CFD) codes and continuously improved through various algorithmic advancements, such as improved shock-capturing techniques \cite{LIU201788,tong2026hybrid}, memory-efficient implicit defect-correction solvers \cite{NakashimaNishikawaLeeCerizza_aiaa_aviation2025}, and gradient-free linear reconstruction schemes \cite{DiskinLiuNishikawa_GradFree_aiaa_aviation2026}. These recent improvements are motivated by the need to achieve efficient computations on graphics processing unit (GPU) platforms and also by the push toward automated CFD simulations with anisotropic viscous grid adaptation \cite{Kleb_etal_aiaa2019-2948,ThompsonNishikawaPadway_aiaa_scitech2023,MoriscoNishikawa_aiaa_scitech2025-0302,Nastac2024ClosedLoop,Liu2025SlidingMesh}. In this paper, we discuss a robustness improvement in the viscous discretization for cell-centered finite-volume methods on unstructured grids.

Viscous schemes play a critical role in algorithmic advancements, particularly for reliable viscous simulations on adaptive unstructured grids. Although they have received less attention than inviscid schemes, robust viscous schemes have often been a concern for unstructured-grid solvers. For example, Ref.~\cite{thomas_diskin_nishikawa:CandF2011} demonstrated that a large edge-derivative contribution in the gradient approximation is essential for multigrid convergence on highly stretched grids. In 2010, a general framework for deriving robust viscous schemes was proposed \cite{nishikawa:AIAA2010}. One particularly successful outcome of this framework is the alpha-damping viscous scheme, which consists of the averaged gradient and a damping term controlled by the damping coefficient $\alpha$. This scheme is a generalized viscous scheme, which includes the classical edge-normal and face-tangent gradient schemes \cite{thomas_diskin_nishikawa:CandF2011,Haselbacher_PhD}, as well as the Mathur-Murthy  scheme \cite{MathurMurthy:NHT1997scheme} as its special cases. Since 2010, the alpha-damping scheme has been found useful especially for practical unstructured-grid CFD codes \cite{nakashima_watanabe_nishikawa:Japan2014,WhiteBaurlePasseSpiegelNishikawa:JANNAF} and continues to attract growing interest for broader deployment \cite{ThompsonNishikawaPadway_aiaa_scitech2023,BATISTIC2026115056,BellostaAbergoNishikawa_aiaa2025-0072}.

Despite its great success, it has been known to practitioners that the alpha-damping scheme sometimes encounters a robustness issue on highly-distorted grids in practical turbulent-flow applications \cite{nakashima_private}. Possible remedies found to be effective by practitioners are to increase the damping coefficient $\alpha$ in the Jacobian and, quite curiously, to change the location of the solution reconstruction at which the solution jump in the damping term is evaluated. Note that the jump is defined by the difference of two solutions linearly reconstructed at an arbitrary nearby location from two adjacent cell centers. In Ref.~\cite{nishikawa:AIAA2010}, two variants of the alpha-damping scheme are proposed for cell-centered finite-volume methods: (1) face-centroid (or face-midpoint in two dimensions) scheme and (2) edge-midpoint scheme. The face-centroid (FC) scheme reconstructs the solutions at a face centroid, where the numerical flux is evaluated, whereas the edge-midpoint (EM) scheme reconstructs them at the midpoint between two adjacent cell centers. Interestingly, it has often been found that switching from the FC scheme to the EM scheme stabilizes implicit solvers  \cite{nakashima_private}, but the mechanism behind this stabilization has not been well understood. In this paper, we will discuss this mechanism in detail, and show that the EM scheme is indeed expected to stabilize the implicit solver for unstructured grids. Although we focus here on second-order cell-centered finite-volume viscous-flux methods, the basic idea potentially has wider implications for inviscid and high-order schemes. 

The main idea discussed here is that if the implicit solver becomes unstable, it may be made stable by reducing the magnitude of the damping term in the residual, or by increasing the magnitude of the damping term from which the Jacobian is derived. For example, increasing the damping coefficient in the Jacobian has the effect of stabilizing the implicit solver. Also, as we will show, a high-order solution reconstruction reduces the damping term in the residual and has a positive impact on stabilizing the implicit solver. Then, we will show that the EM scheme has the effect of reducing the magnitude of the damping term on unstructured grids, and therefore, improves convergence characteristics for unstructured-grid calculations over the FC scheme. Ref.~\cite{nishikawa_nakashima_watanabe:jcp2017} presents a similar study, but it solely examines the effect of the damping coefficient, not paying attention to the relative magnitude of the damping terms in the residual and Jacobian. In this paper, we extend the study and propose to focus on the relative magnitude of the damping terms. 

The next section describes the target cell-centered finite-volume method and the alpha-damping viscous scheme. Then, in Section 3, we demonstrate by a Fourier analysis that an implicit defect-correction solver can become unstable if the damping term in the residual is too large, or equivalently if the lower-order jump used to derive the Jacobian becomes too small. In Section 4, we demonstrate that the damping term computed with the EM scheme can be smaller in magnitude than that computed with the FC scheme. In Section 5, we present numerical results using unstructured grids that support the analysis and demonstrate the effectiveness of the EM scheme. Finally, in Section 6, we conclude the paper with remarks.

\section{Second-Order Finite-Volume Solver}
\label{intro} 

\indent

\begin{figure}[t]
\begin{center}
\begin{minipage}[b]{0.5\textwidth}
\begin{center}
\includegraphics[width=0.99\textwidth,trim=0 0 0 0,clip]{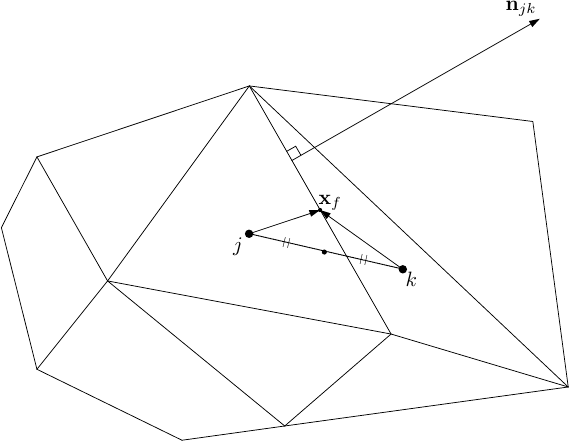}
\caption{A mixed-element grid for a cell-centered finite-volume discretization. }
\label{fig:ccfv_res_stencil}
\end{center}
\end{minipage}
\end{center}
\end{figure}

Consider the steady compressible Navier-Stokes equations:
\begin{eqnarray}
 \mbox{div} \left[ {\cal F}^{inv}({\bf w}) +  {\cal F}^{vis}({\bf w}, \nabla {\bf w} ) \right] = {\bf 0},
\label{diff_form}
\end{eqnarray}
where ${\cal F}^{inv}$ and ${\cal F}^{vis}$ are the inviscid and viscous flux tensors, and ${\bf w}$ is a vector of primitive variables. Its second-order cell-centered finite-volume discretization at a cell $j$ on an unstructured mixed-element grid (see Figure \ref{fig:ccfv_res_stencil}) is given by
\begin{eqnarray}
  \frac{1}{V_j} \sum_{k \in \{ k_j\} } 
  \left[ {\Phi}^{inv}_{jk} + {\Phi}^{vis}_{jk} \right]
  |{\bf n}_{jk}| = {\bf 0}, 
\label{fvm_diff}
\end{eqnarray}
where $V_j$ is the volume of the cell $j$, $\{ k_j\}$ is a set of face neighbors of the cell $j$, ${\bf n}_{jk}$ is the scaled face-normal vector, and ${\Phi}^{inv}_{jk}$ and ${\Phi}^{vis}_{jk}$ are numerical fluxes for the inviscid and viscous fluxes, respectively,  projected along the face-normal direction. The inviscid flux ${\Phi}^{inv}_{jk} = {\Phi}^{inv}_{jk}({\bf w}_L,{\bf w}_R)$ is computed at the face centroid based on two states, ${\bf w}_L$ and ${\bf w}_R$, linearly reconstructed at the face centroid ${\bf x}_f$ from the centroids of the cells $j$ and $k$, respectively,
\begin{eqnarray}
{\bf w}_L = {\bf w}_j + \overline{\nabla} {\bf w}_j \cdot 
(  {\bf x}_f - {\bf x}_j ) ,
\quad
{\bf w}_R = {\bf w}_k + \overline{\nabla} {\bf w}_k \cdot 
(  {\bf x}_f - {\bf x}_k ),
\end{eqnarray}
where $\overline{\nabla} {\bf w}_j$ and $\overline{\nabla} {\bf w}_k$ are gradients at the cells $j$ and $k$, respectively, computed by a linear least-squares (LSQ) method {\color{black} (see, e.g., Ref.\cite{nishikawa_stencil:JCP2019})}. In this paper, the inviscid flux is the robust Roe flux \cite{Nishikawa_RobustFluxes:jcp2020}. 

The viscous flux ${\Phi}^{vis}_{jk} = {\Phi}^{vis}_{jk}({\bf w}_j,{\bf w}_k, \left. \nabla {\bf w} \right|_f)$ is computed based on the face gradient $\left. \nabla {\bf w} \right|_f$ evaluated by the alpha-damping scheme \cite{nishikawa:AIAA2010}, which is the main subject of this paper: e.g., for the $x$-component of the velocity, $u$, it is given by
\begin{eqnarray}
\left. \nabla u \right|_f
 = 
 \frac{ \overline{\nabla} u_j + \overline{\nabla} u_k }{2} 
 +
 \frac{\alpha}{ | ( {\bf x}_k - {\bf x}_j )
 \cdot \hat{\bf n}_{jk} |} 
 ( u_R - u_L ) \,  \hat{\bf n}_{jk} ,
 \label{alpha_damp_u}
\end{eqnarray}
where $ \hat{\bf n}_{jk} = {\bf n}_{jk} / |{\bf n}_{jk}|$. Here, the second term is called the damping term, serving to damp high-frequency errors (see Ref.~\cite{nishikawa:AIAA2010}), and $\alpha$ is called the damping coefficient. Throughout the paper, we set $\alpha = 4/3$ unless otherwise stated, which gives fourth-order accuracy for diffusion in one dimension \cite{nishikawa:AIAA2010} and has also been found accurate even on unstructured grids \cite{jalali_etal:CF2014}.
The same gradient scheme is applied to the other velocity components and the temperature required in the viscous flux. The viscosity is defined by Sutherland's law and evaluated by the temperature given by the arithmetic average of the cell center values, which is sufficient for second-order accuracy (see Ref.~\cite{NishikawaDiskin_Note_ViscosityAveraging:2022}).

Our interest is in the stability of an implicit defect-correction solver, widely employed in practical CFD codes:
\begin{eqnarray}
{\bf U}^{n+1} = {\bf U}^{n} + \Delta {\bf U},
\end{eqnarray}
where ${\bf U}$ is a global vector of numerical solutions in the conservative variables, $n$ is the iteration counter, and $ \Delta {\bf U}$ is computed by relaxing the linearized system:
\begin{eqnarray}
\frac{ \partial \overline{\bf Res} }{\partial   {\bf U} }  \Delta {\bf U} = -{\bf Res} ( {\bf U}^n),
\end{eqnarray}
where $\overline{\bf Res}$ is the global vector of the residual where all LSQ gradients are ignored (i.e., no reconstruction and no gradients in the viscous scheme). The linear system is relaxed by a multi-color Gauss-Seidel relaxation scheme. For the viscous part, which is our focus here, the residual Jacobian is constructed by differentiating only the damping terms proportional to the jump between two cell-center solutions.  As typically done, a diagonal matrix of pseudotime steps, controlled by the CFL number, is added to the Jacobian for robustness. Our focus is on the impact of the damping term on the stability of the implicit defect-correction solver.


\section{Analysis in One Dimension}
\label{intro}

To gain insight into the stability of the implicit defect-correction solver, we analyze it for diffusion in one dimension, $ \partial_{xx} u = f(x)$. Consider a uniform grid of spacing $h$ in a unit domain with cell centers defined by $x_j = (j-1/2) h $, where $j=1,2,3, \cdots, N$, and $h = 1/N$. For the purpose of analysis, we ignore boundary effects by specifying exact solution values at cells near the left and right boundaries (i.e., $j=1, 2, N-1$, and $N$). The residual at cell $j$ is given by
\begin{eqnarray}
 Res_j 
=
  f_{j+1/2} - f_{j-1/2}  - f(x_j) h,
\label{oned_fv_res_j}
\end{eqnarray}
where 
\begin{eqnarray}
f_{j+1/2} 
=
 \frac{ ( \partial_x u)_j + ( \partial_x u)_{j+1}  }{2} 
 +
 \frac{\alpha}{ h } 
 ( u_R - u_L ), 
 \quad
 u_L = u_j     + \frac{h}{2} ( \partial_x u)_j , 
 \quad
 u_R = u_{j+1} - \frac{h}{2} ( \partial_x u)_{j+1} , 
\label{oned_fv_alpha_flux}
\end{eqnarray}
\begin{eqnarray}
 ( \partial_x u)_j = \frac{ u_{j+1} - u_{j-1}}{2h},
 \quad
 ( \partial_x u)_{j+1} = \frac{ u_{j+2} - u_{j}}{2h}.
 \label{cd_grad}
\end{eqnarray}
The implicit defect-correction solver is given by
\begin{eqnarray}
u_j^{n+1} = u_j^n +  \Delta u_j ,
\label{oned_implicit_solver_update}
\end{eqnarray}
for $j=3,4, \cdots, N-2$, where $n=0,1,2,\cdots$ is the iteration counter, where the linearized system of the implicit defect-correction solver is given at a cell $j$ by
\begin{eqnarray}
\frac{\alpha_{Jac}}{h} 
\Delta u_{j-1}
-
\frac{2 \alpha_{Jac}}{h} 
\Delta u_j 
+
\frac{\alpha_{Jac}}{h} 
\Delta u_{j+1}
=
- Res_j,
\label{oned_implicit_system}
\end{eqnarray}
for $j=3,4, \cdots, N-2$. Note that we have introduced $\alpha_{Jac}$ to allow us to use a different value of the damping coefficient in the Jacobian. In this one-dimensional solver, the linear system is inverted directly. Furthermore, the pseudotime terms are omitted, which corresponds to an infinite CFL number.

To analyze the stability of the implicit defect-correction solver, we consider a Fourier mode, $U_0 \exp( x \beta /h )$, where $U_0$ is an amplitude and $\beta$ is a frequency. Inserting this into Equation (\ref{oned_implicit_solver_update}) with Equation (\ref{oned_implicit_system}), we obtain 
\begin{eqnarray}
U_0^{n+1} 
=
s(\beta) 
U_0^{n} ,
\quad
s(\beta)  = 1  - \frac{\lambda_{RHS}}{\lambda_{LHS}},
\label{s_beta}
\end{eqnarray}
where $\lambda_{LHS}$ and $\lambda_{RHS}$ are Fourier symbols of the Jacobian and residual operators, respectively:
\begin{eqnarray}
\lambda_{LHS} 
=
 \frac{ 2 \alpha_{Jac}  (  \cos \beta - 1  ) }
 { h }, 
 \quad
\lambda_{RHS} 
=
 \frac{\alpha  ( \sin^2 \beta  + 2 \cos \beta - 2  ) - \sin^2 \beta  }
 { h }.
 \label{fourtier_symbols}
\end{eqnarray}
The solver is stable if 
\begin{eqnarray}
s_{max} = \max_{ \beta \in [ 0, \pi ] } | s(\beta) | < 1 .
\end{eqnarray}
It is important to note in Equation (\ref{s_beta}) that if $|\lambda_{RHS}| < |\lambda_{LHS}|$, then the solver is stable. This implies that the solver, if unstable, may be stabilized by reducing the damping term in the residual or by increasing the damping coefficient in the Jacobian. 

\begin{figure}[htbp!]
    \centering
          \begin{subfigure}[t]{0.49\textwidth}
    \centering
        \includegraphics[width=0.99\textwidth,trim=0 0 0 0,clip]{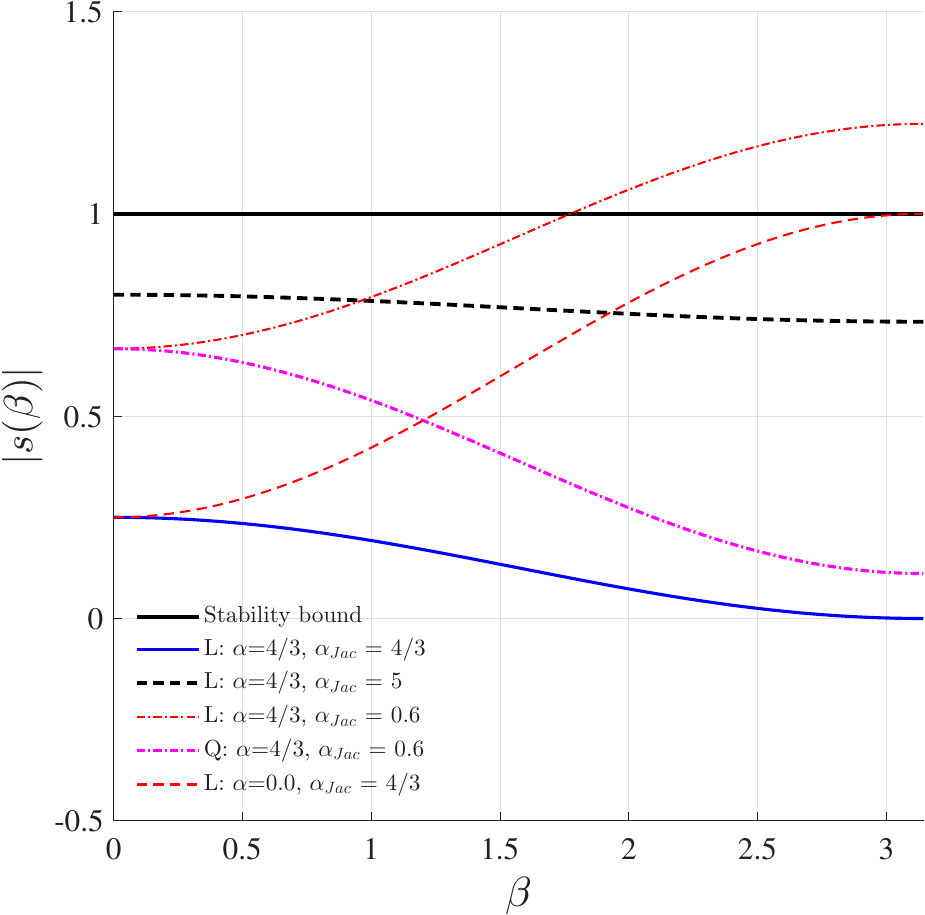}
          \caption{Stability analysis results.}
          \label{fig:oned_analysis}
      \end{subfigure}
      \hfill
          \begin{subfigure}[t]{0.49\textwidth}
    \centering
        \includegraphics[width=0.99\textwidth,trim=0 0 0 0,clip]{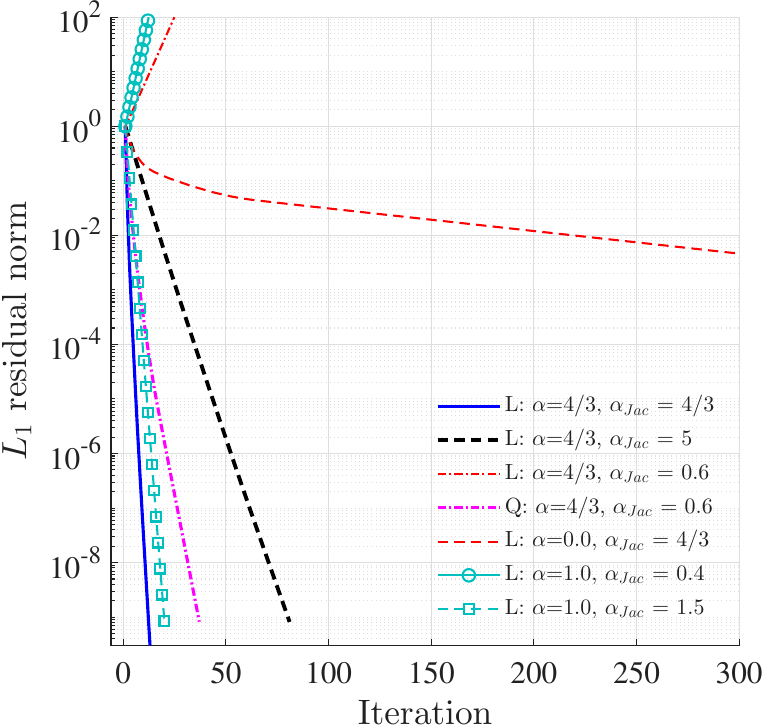}
          \caption{Iterative convergence.}
          \label{fig:oned_res_conv}
      \end{subfigure}
            \caption{
          \label{edge_collapse_idea_full}
Results for one-dimensional problem. 
} 
\end{figure}

The function $s(\beta)$ is plotted in Figure \ref{fig:oned_analysis}, where L indicates that the gradient is computed by the central difference formula (\ref{cd_grad}). For $\alpha=4/3$, it is observed, as expected, that the solver is stable with $\alpha_{Jac}=4/3$ and also with a larger value ($\alpha_{Jac}=5$) although convergence is expected to be slower. We also observe that the solver becomes unstable with $\alpha_{Jac}=0.6$. In fact, noticing that $s_{max}$ reaches 1 at $\beta=\pi$ with decreasing $\alpha_{Jac}$, we can derive that the solver becomes unstable when $\alpha_{Jac} \le \alpha /2 = 2/3 = 0.666 \cdots$, for $\alpha=4/3$. 

%
%

As pointed out earlier, we expect the implicit solver to become more stable when the damping term in the residual is smaller in magnitude. This can be achieved by a high-order solution reconstruction. To demonstrate this, we replace the linear reconstruction by the following quadratic reconstruction:
\begin{eqnarray}
 u_L = u_j     + \frac{h}{2} ( \partial_x u)_j 
+ \frac{1}{2} ( \partial_{xx} u)_j \left( \frac{h}{2} \right)^2 , 
 \quad
 u_R = u_{j+1} - \frac{h}{2} ( \partial_x u)_{j+1} 
 + \frac{1}{2} ( \partial_{xx} u)_{j+1}  \left( \frac{h}{2} \right)^2 , 
\label{oned_fv_alpha_flux_quadratic}
\end{eqnarray}
where the second derivatives are computed by the three-point central difference formula for simplicity. For a smooth function, it can be shown that the quadratic reconstruction gives
 \begin{eqnarray}
 u_R - u_L = - \frac{1}{8} \left( \partial_{xxx} u_{j+1/2}  \right)h^3 + O(h^5),
\end{eqnarray}
while the linear reconstruction gives
 \begin{eqnarray}
 u_R - u_L = - \frac{1}{4} \left( \partial_{xxx} u_{j+1/2} \right) h^3 + O(h^5).
\end{eqnarray}
Therefore, the quadratic reconstruction yields a smaller jump. The quadratic reconstruction gives a slightly different Fourier symbol for the residual:
\begin{eqnarray}
\lambda_{RHS} 
=
 \frac{\alpha  ( \sin^2 \beta  + 2 \cos \beta - 2  ) - 2 \sin^2 \beta  }
 { h },
\end{eqnarray}
which differs from the original $\lambda_{RHS}$ in Equation (\ref{fourtier_symbols}) by a factor of two in front of $\sin^2 \beta$. In Figure \ref{fig:oned_analysis}, this case is indicated by Q. As can be seen, the case of $\alpha_{Jac}=0.6$, which is unstable with the linear reconstruction, turns stable with the quadratic  reconstruction. However, as one can expect, the damping term is essential to stability and therefore, the solver cannot be stable with a too small damping term. To show this, we added one more data set corresponding to $\alpha=0$, i.e., zero damping term in the residual. As can be seen, $s_{max}=1$ at $\beta=\pi$ and indeed, the solver is expected to converge extremely slowly if it ever converges. 

To confirm the analysis, we actually solve the residual equations for all the cases considered in the analysis, on a uniform grid with $N=16$. The exact solution is given by $u_{exact}(x) = \sin(2 \pi x) + x^2/2$, and the forcing term is set to $f(x) = \partial_{xx} u_{exact}(x)$. Results are shown in Figure \ref{fig:oned_res_conv}. As predicted by the analysis, the implicit solver converged all the cases except for the linear reconstruction case with $\alpha=4/3$ and $\alpha=0.6$. Also, as expected, the solver seems to be converging, but extremely slowly in the zero-damping-term case ($\alpha=0$). 

Finally, we consider a case, where the damping term is evaluated from solutions at cell centers, i.e., $u_R-u_L = u_{j+1} - u_j$, in the residual. Note that this solution jump is called the unreconstructed jump in the rest of the paper. Then, with $\alpha=1$, the residual reduces to the three-point finite-difference scheme for diffusion, $( u_{j+1} - 2 u_j + u_{j-1} )/h^2$. 
In this case, we obtain
\begin{eqnarray}
| s(\beta) | 
=  
\left| 
1  - \frac{\alpha}{\alpha_{Jac}}
\right|.
\end{eqnarray}
Assuming that both damping coefficients are positive, as required for the damping term to damp high-frequency errors, we see that the implicit solver is stable if 
\begin{eqnarray}
\alpha_{Jac} > \frac{\alpha}{2} = \frac{1}{2},
\end{eqnarray}
where $\alpha=1$ is assumed because that is the only choice that makes the residual consistent and second-order accurate. This implies that the solver may be stabilized by reverting to the cell-center-solution jump with a sufficiently large $\alpha_{Jac}$. Two examples are shown in Figure \ref{fig:oned_res_conv}: $\alpha=1$ with $\alpha_{Jac} = 0.4$ and $\alpha_{Jac} = 1.5$. As expected, the implicit solver diverges in the former but converges in the latter. In one dimension, the residual remains second-order accurate with $\alpha=1$ but the alpha-damping scheme with the unreconstructed jump is not even first order accurate in two and three dimensions \cite{Nishikawa_RobustFluxes:jcp2020}. Nevertheless, it might be a useful option to have in cases that one wishes to obtain a numerical solution even if it is not as accurate as one would hope for. Later, we will show some examples.

The one-dimensional study indicates that the implicit solver may be stabilized by increasing $\alpha$ in the Jacobian or by making the solution jump smaller (e.g., by a high-order solution reconstruction). With this in mind, we return to the two-dimensional unstructured-grid solver and discuss robust implicit solver convergence in conjunction with the FC and EM schemes.

\section{Face-Centroid and Edge-Midpoint Schemes}
\label{intro}

\indent

In two dimensions, we focus on the damping term and its impact on solver stability. Note that the damping term does not affect the order of accuracy as long as it is designed to vanish for a linear function \cite{nishikawa:AIAA2010}. Therefore, it may be evaluated at a location other than the face centroid. Then, we seek a position that potentially reduces the magnitude of the damping term. For our purpose, it suffices to consider the solution jump:
\begin{eqnarray}
  \Delta u  \equiv  u_R - u_L,
   \quad
   u_R =  u_k  +  \left(  {\bf x}_{a} -   {\bf x}_{k}  \right)   \cdot \overline{\nabla} u_k,
  \quad
   u_L = u_j   +  \left(  {\bf x}_{a} -   {\bf x}_{j}  \right)  \cdot \overline{\nabla} u_j ,
   \label{u_edgemidpoint}
\end{eqnarray}
where $u$ is a solution variable and ${\bf x}_{a}$ is the position of the reconstruction. Assuming that the gradients are first-order accurate, which is typical in linear LSQ gradients on unstructured grids, we expand the solutions and gradients:
\begin{eqnarray}
   u_k
   &=&
   u_a 
   +  \left(  {\bf x}_{k} -   {\bf x}_{a}  \right)  \cdot \nabla u_a 
   + \frac{1}{2}  \left[ \left(  {\bf x}_{k} -   {\bf x}_{a}  \right)  \cdot \nabla \right]^2 u_a     + O(h^3) 
    \\  [2ex] 
   u_j
   &=&
   u_a 
   +  \left(  {\bf x}_{j} -   {\bf x}_{a}  \right)  \cdot \nabla u_a 
   + \frac{1}{2}  \left[ \left(  {\bf x}_{j} -   {\bf x}_{a}  \right)  \cdot \nabla \right]^2 u_a     + O(h^3) ,
\end{eqnarray}
\begin{eqnarray}
\overline{\nabla}  u_k
   &=&
 \nabla  u_a
   +   \left[   \left(  {\bf x}_{k} -   {\bf x}_{a}  \right)  \cdot \nabla \right] \nabla u_a  + {\bf C}_k h+ O(h^2) ,
   \\  [2ex]    
\overline{\nabla}  u_j
   &=&
 \nabla  u_a
   +   \left[   \left(  {\bf x}_{j} -   {\bf x}_{a}  \right)  \cdot \nabla \right] \nabla u_a + {\bf C}_j h  + O(h^2) ,
\end{eqnarray}
where ${\bf C}_j h$ and ${\bf C}_k h$ are first-order gradient errors at $j$ and $k$, respectively.  Then, we find
\begin{eqnarray}
  \Delta u 
&=&
   \left[ \left(  \frac{ {\bf x}_{j} + {\bf x}_{k} }{2}  - {\bf x}_{a}  \right)  \cdot \nabla \right] 
   \left[ \left(  {\bf x}_{k} -   {\bf x}_{j}  \right)  \cdot \nabla \right] 
   u_a    
   + 
   \left[   \left(  {\bf x}_{k} -   {\bf x}_{a}  \right)  \cdot \nabla \right] 
   \left[   \left(  {\bf x}_{a} -   {\bf x}_{k}  \right)  \cdot \nabla \right] 
   u_a  
   \\ [2ex]
&-&
   \left[   \left(  {\bf x}_{j} -   {\bf x}_{a}  \right)  \cdot \nabla \right] 
   \left[   \left(  {\bf x}_{a} -   {\bf x}_{j}  \right)  \cdot \nabla \right] 
   u_a     
   +  \left(  {\bf x}_{a} -   {\bf x}_{k}  \right)  \cdot  {\bf C}_k h 
   -  \left(  {\bf x}_{a} -   {\bf x}_{j}  \right)  \cdot  {\bf C}_j h 
   + O(h^3)
    \\ [2ex]
&=& 
- \left[ \left(  \frac{ {\bf x}_{j} + {\bf x}_{k} }{2}  - {\bf x}_{a}  \right)  \cdot \nabla \right] 
   \left[ \left(  {\bf x}_{k} -   {\bf x}_{j}  \right)  \cdot \nabla \right] 
   u_a   
   +  \left(  {\bf x}_{a} -   {\bf x}_{k}  \right)  \cdot  {\bf C}_k h 
   -  \left(  {\bf x}_{a} -   {\bf x}_{j}  \right)  \cdot  {\bf C}_j h 
   + O(h^3),
\end{eqnarray}
and thus 
\begin{eqnarray}
| \Delta u |
  \le 
\left| 
\left[ \left(  \frac{ {\bf x}_{j} + {\bf x}_{k} }{2}  - {\bf x}_{a}  \right)  \cdot \nabla \right] 
   \left[ \left(  {\bf x}_{k} -   {\bf x}_{j}  \right)  \cdot \nabla \right] 
   u_a   
\right|
+
\left| 
   \left(  {\bf x}_{a} -   {\bf x}_{k}  \right)  \cdot  {\bf C}_k h 
\right|
+
\left| 
  \left(  {\bf x}_{a} -   {\bf x}_{j}  \right)  \cdot  {\bf C}_j h 
\right| + \cdots.
\label{du_bound}
\end{eqnarray}
Note that all three terms on the right hand side are of $O(h^2)$, and clearly, the first term of Equation (\ref{du_bound}) can be eliminated by choosing the reconstruction position to be the edge midpoint:
 \begin{eqnarray}
   {\bf x}_{a}  =  \frac{  {\bf x}_{j}  +  {\bf x}_{k}  }{2}. 
\end{eqnarray}
This shows that the damping term has a smaller upper bound for the EM scheme than the FC scheme. 
This appears to explain why the EM scheme has been found more robust than the FC scheme in practical unstructured-grid simulations. 

Another observation is that the upper bound can be reduced by using a quadratic LSQ gradient method, which eliminates the second and third terms on the right-hand side of Equation (\ref{du_bound}). It is also possible, as suggested by the one-dimensional analysis and will be demonstrated later, that the implicit solver is stabilized by increasing the damping coefficient in the Jacobian.


Finally, the main point of the present study is that the implicit solver tends to be more stable if the damping term in the residual is forced to be smaller relative to the term in the Jacobian, and this can be achieved, for example, by the techniques mentioned earlier. In particular, the study has provided a reason that implicit solvers are expected to be more robust with the EM scheme, and it deserves to be a default choice.

\section{Results}
\label{results} 

\indent

\subsection{Laminar Flow over a Flat Plate}
\label{results_laminar_fp} 

%
   \begin{figure}[htbp!]
    \centering
          \begin{subfigure}[t]{0.49\textwidth}
    \centering
        \includegraphics[width=0.99\textwidth,trim=0 0 0 0,clip]{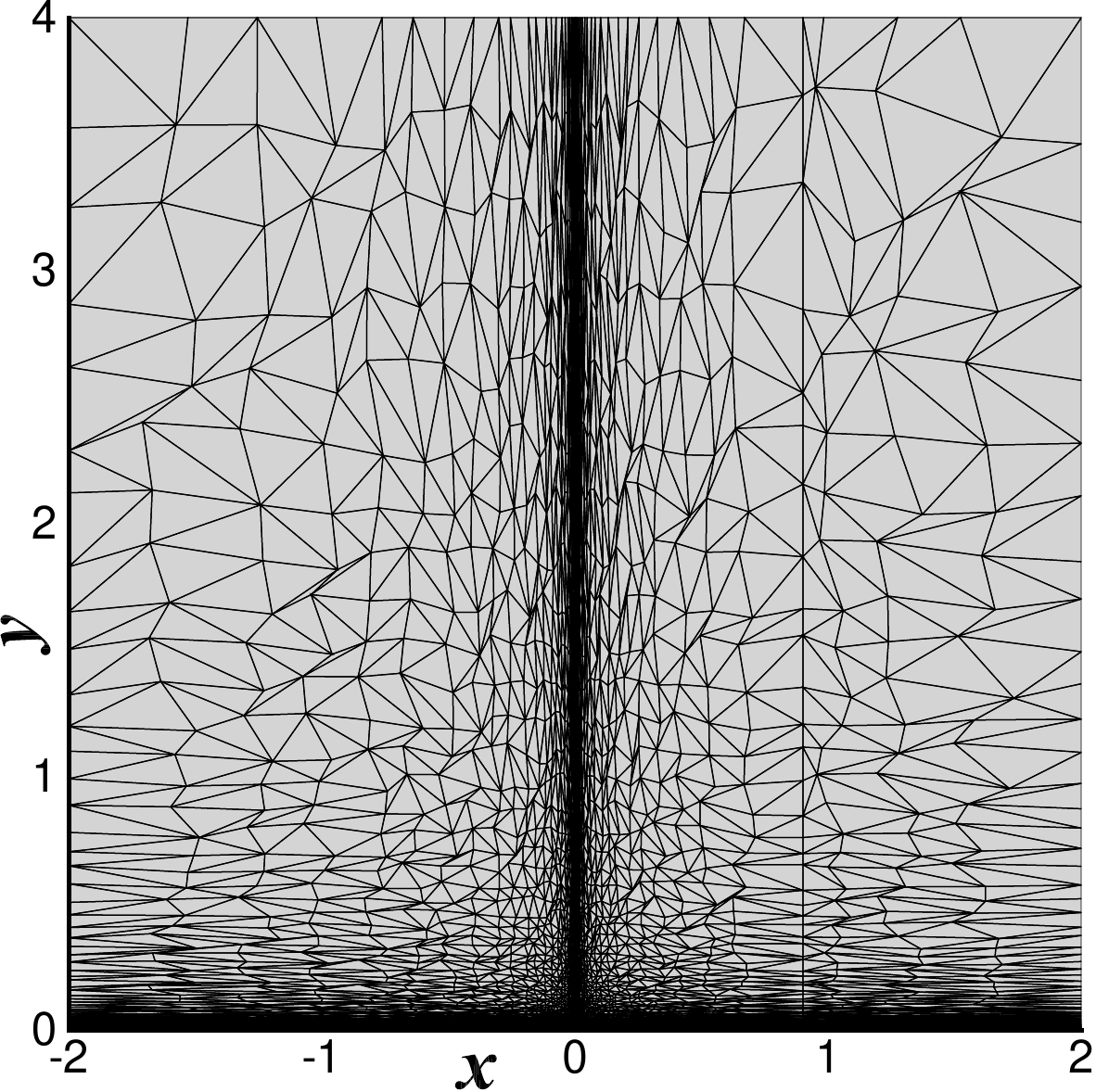}
          \caption{Mixed-element grid.}
          \label{fig:twod_fp_grid}
      \end{subfigure}
      \hfill
          \begin{subfigure}[t]{0.49\textwidth}
    \centering
        \includegraphics[width=0.99\textwidth,trim=0 0 0 0,clip]{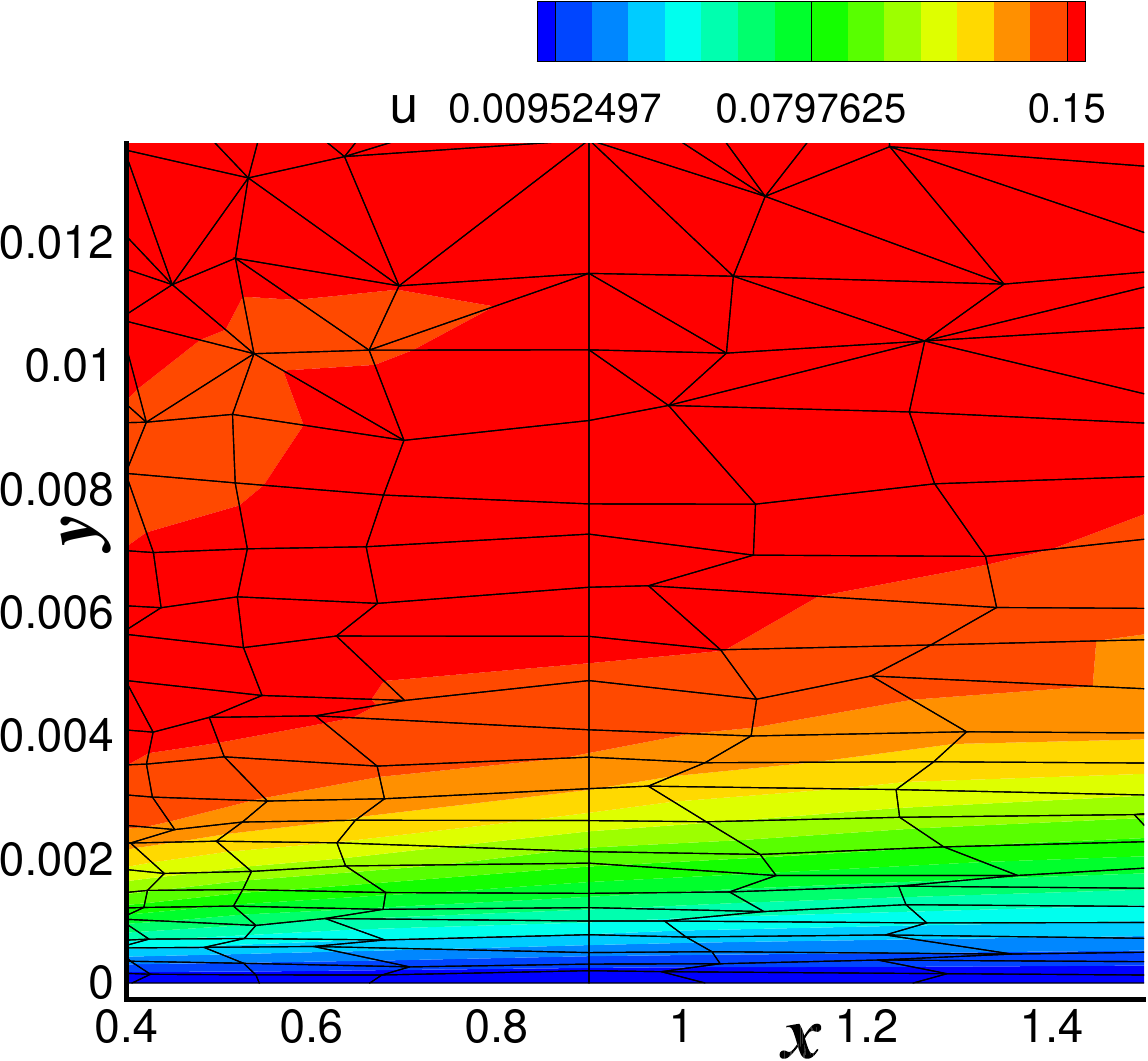}
          \caption{Zoomed-in view.}
          \label{fig:twod_fp_grid_zoomed}
      \end{subfigure}
            \caption{
          \label{edge_collapse_idea_full}
The grid used in the flat plate case.
} 
\end{figure}
%

We consider a laminar flow over a flat plate at $M_\infty = 0.15$ and $Re_\infty = 10^6$ in the domain $\{(x,y) |   -2 \le x \le 2, \, 0 \le y \le 4 \}$. The bottom right half boundary, $ 0 \le x \le 2$, is taken to be a flat plate, and the no-slip wall condition is applied weakly through the right state to the numerical flux at each boundary face. Namely, the velocity in the right state is defined by the velocity reconstructed from inside with a negative sign, so that the velocity is zero on average. The pressure is copied from the interior and the density is defined by $\gamma p_L / T_{wall}$, where $\gamma=1.4$ and $T_{wall}=300[K]$. Other boundary conditions are also enforced weakly (see, e.g., Ref.~\cite{liu_nishikawa_aiaa2016-3969}): a symmetry condition on the bottom left boundary; a free stream condition on the left boundary, a far-field condition on the top boundary; a back-pressure outflow condition on the right boundary face. The grid is an irregular mixed-element grid with 3,525 nodes, 5,048 triangles, and 880 quadrilaterals, as shown in Figure \ref{fig:twod_fp_grid}. See also a zoomed-in view in Figure \ref{fig:twod_fp_grid_zoomed}, where the $x$-velocity contours are shown and a transition from quadrilateral cells to triangular cells can be seen. Note that there is a straight vertical grid line at $x=0.9$, which is introduced for a solution sampling purpose.

In the implicit defect-correction solver, we perform the linear relaxation until the linear-system residual is reduced by one order of magnitude or until the maximum of 150 relaxations is reached. The initial values are set by a free stream condition with random perturbation introduced to all variables in all cells. The CFL number for the pseudotime term is set to $10^{15}$, and it will be reduced only if the residual norm increases more than one order of magnitude in one iteration or unrealizable states are encountered in the solution update. The solver is taken to be converged when the maximum $L_1$ norm of the residuals is reduced by six orders of magnitude. 

\begin{figure}[htbp!]
    \centering
          \begin{subfigure}[t]{0.32\textwidth}
    \centering
        \includegraphics[width=0.99\textwidth,trim=0 0 0 0,clip]{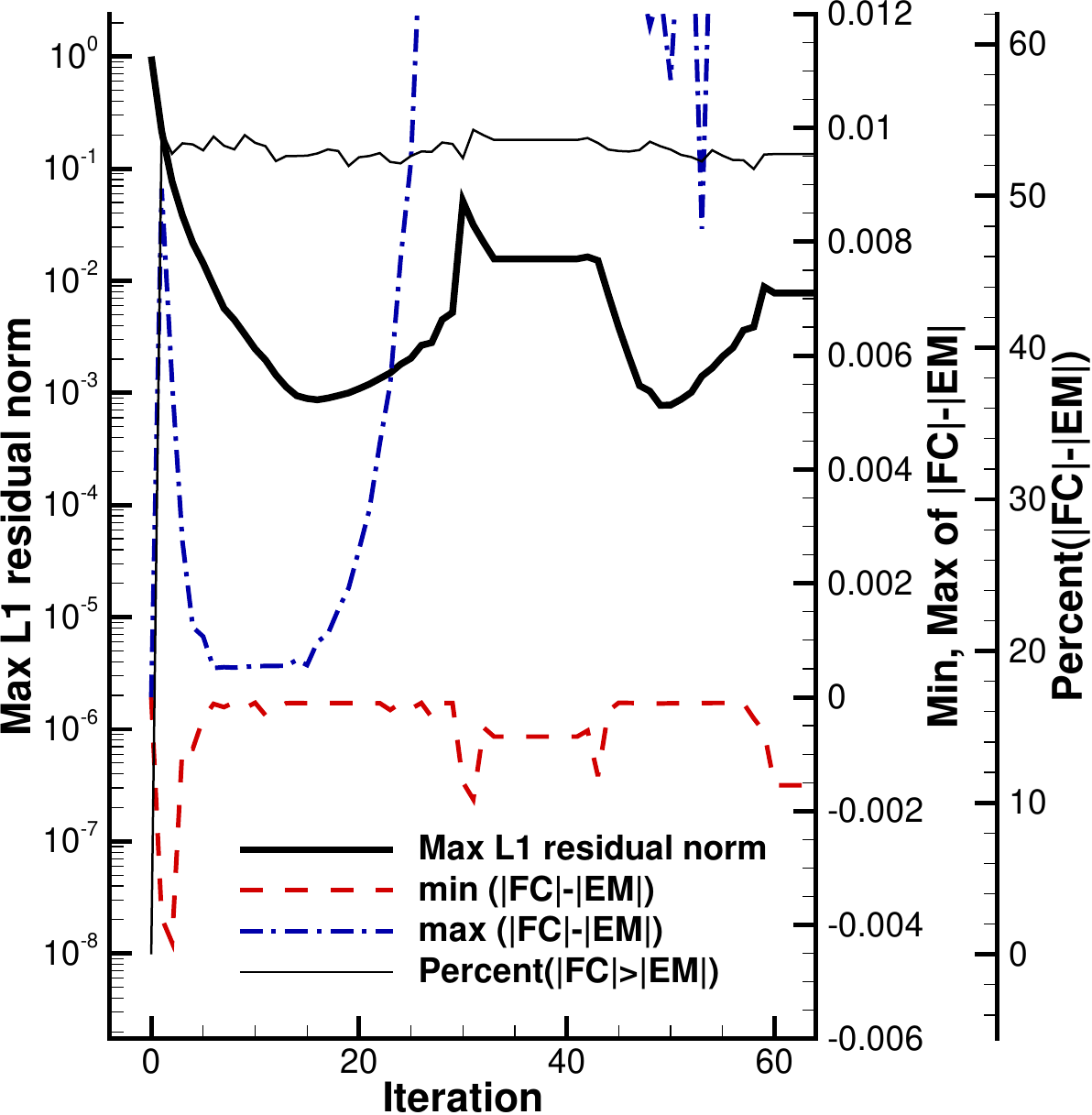}
          \caption{FC: $\alpha=4/3$, $\alpha_{Jac}=4/3$.}
          \label{fig:twod_conv_fc_4o3}
      \end{subfigure}
      \hfill
          \begin{subfigure}[t]{0.32\textwidth}
    \centering
        \includegraphics[width=0.99\textwidth,trim=0 0 0 0,clip]{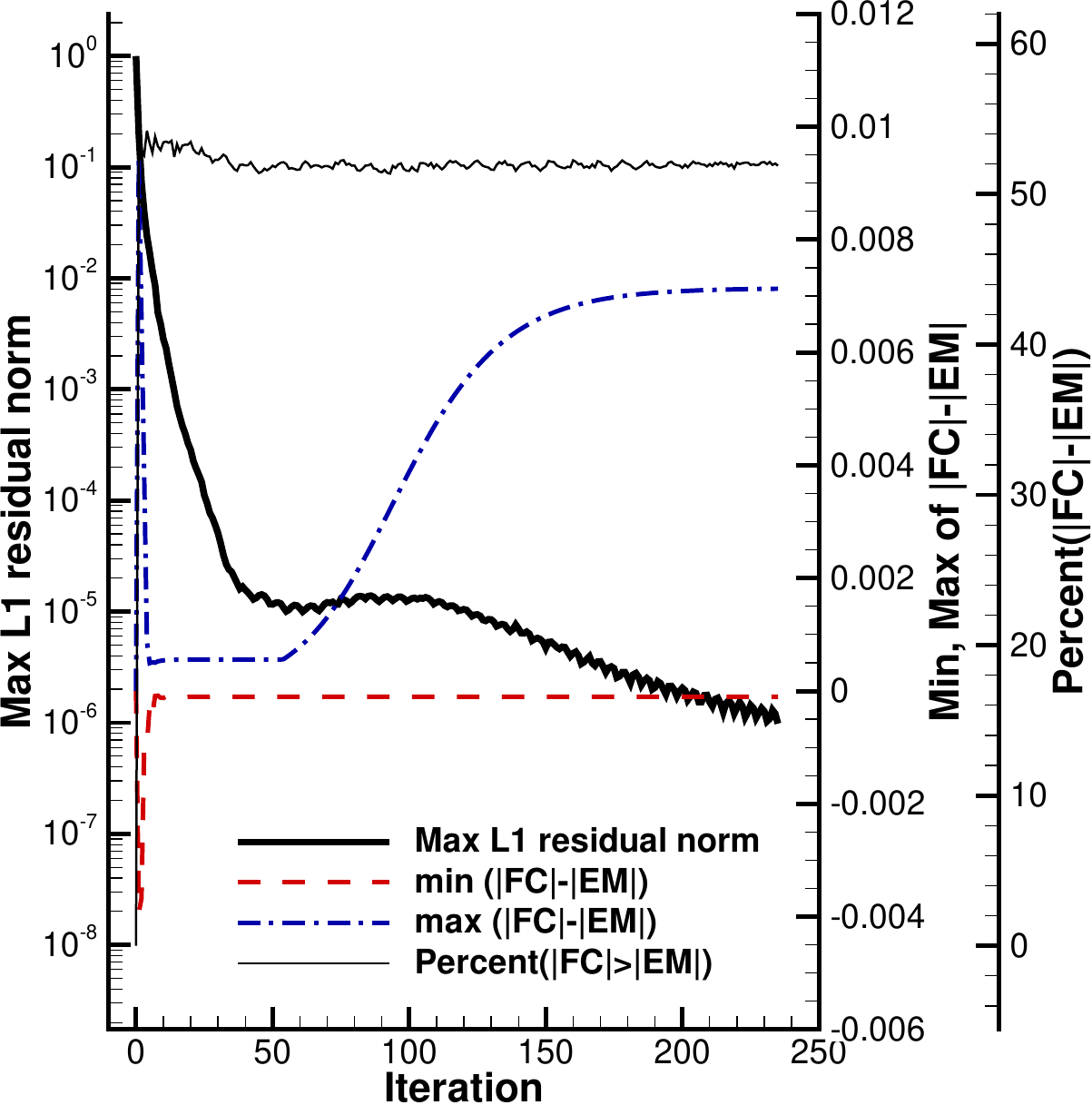}
          \caption{FC: $\alpha=4/3$, $\alpha_{Jac}=5$.}
          \label{fig:twod_conv_fc_4o3_largeJac}
      \end{subfigure}
      \hfill
          \begin{subfigure}[t]{0.32\textwidth}
    \centering
        \includegraphics[width=0.99\textwidth,trim=0 0 0 0,clip]{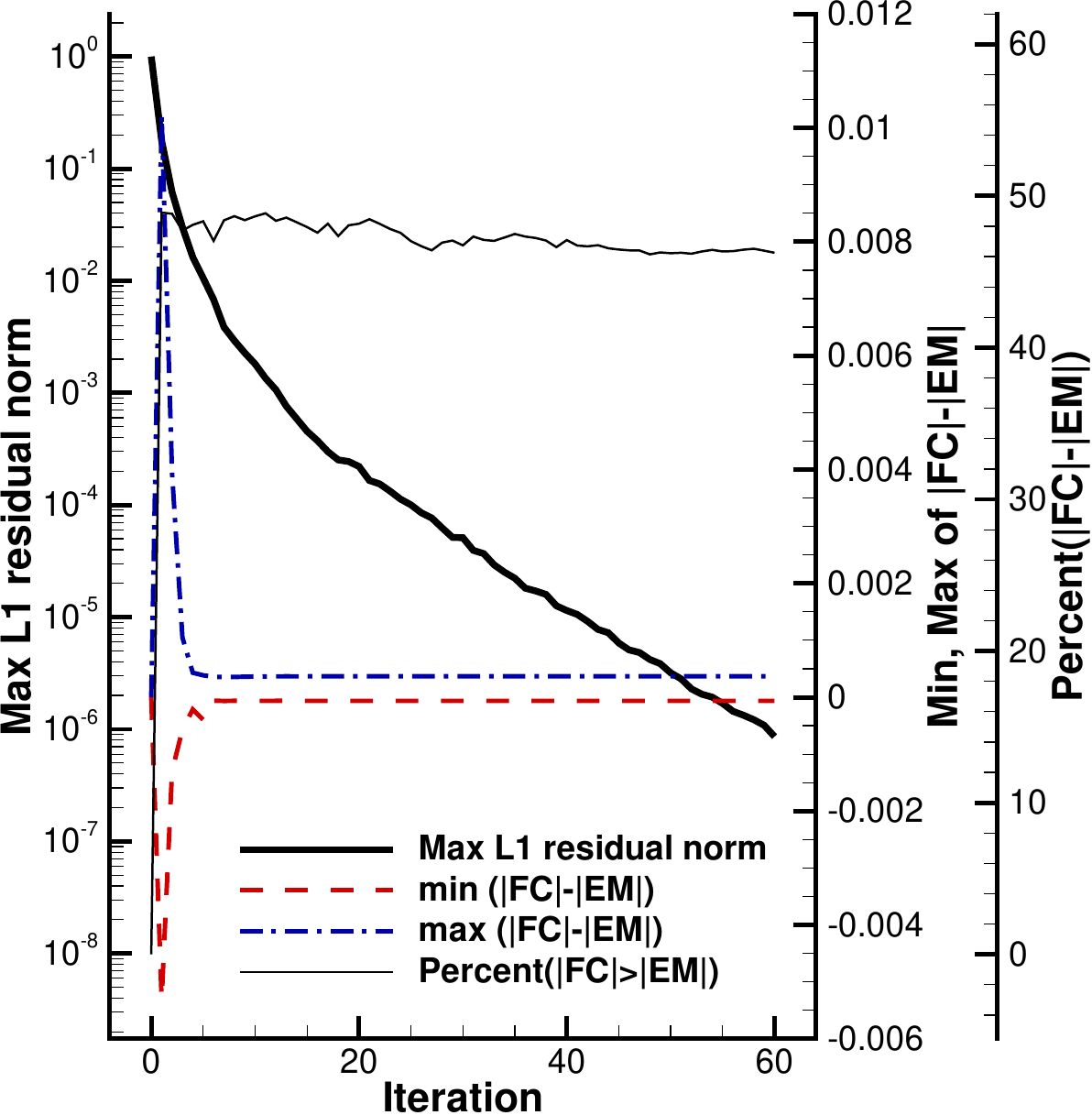}
          \caption{FC (QLSQ): $\alpha=4/3$, $\alpha_{Jac}=4/3$.}
          \label{fig:twod_conv_fc_4o3_qlsq}
      \end{subfigure}
      \\ \vspace{0.3cm}
    \centering
          \begin{subfigure}[t]{0.32\textwidth}
    \centering
        \includegraphics[width=0.99\textwidth,trim=0 0 0 0,clip]{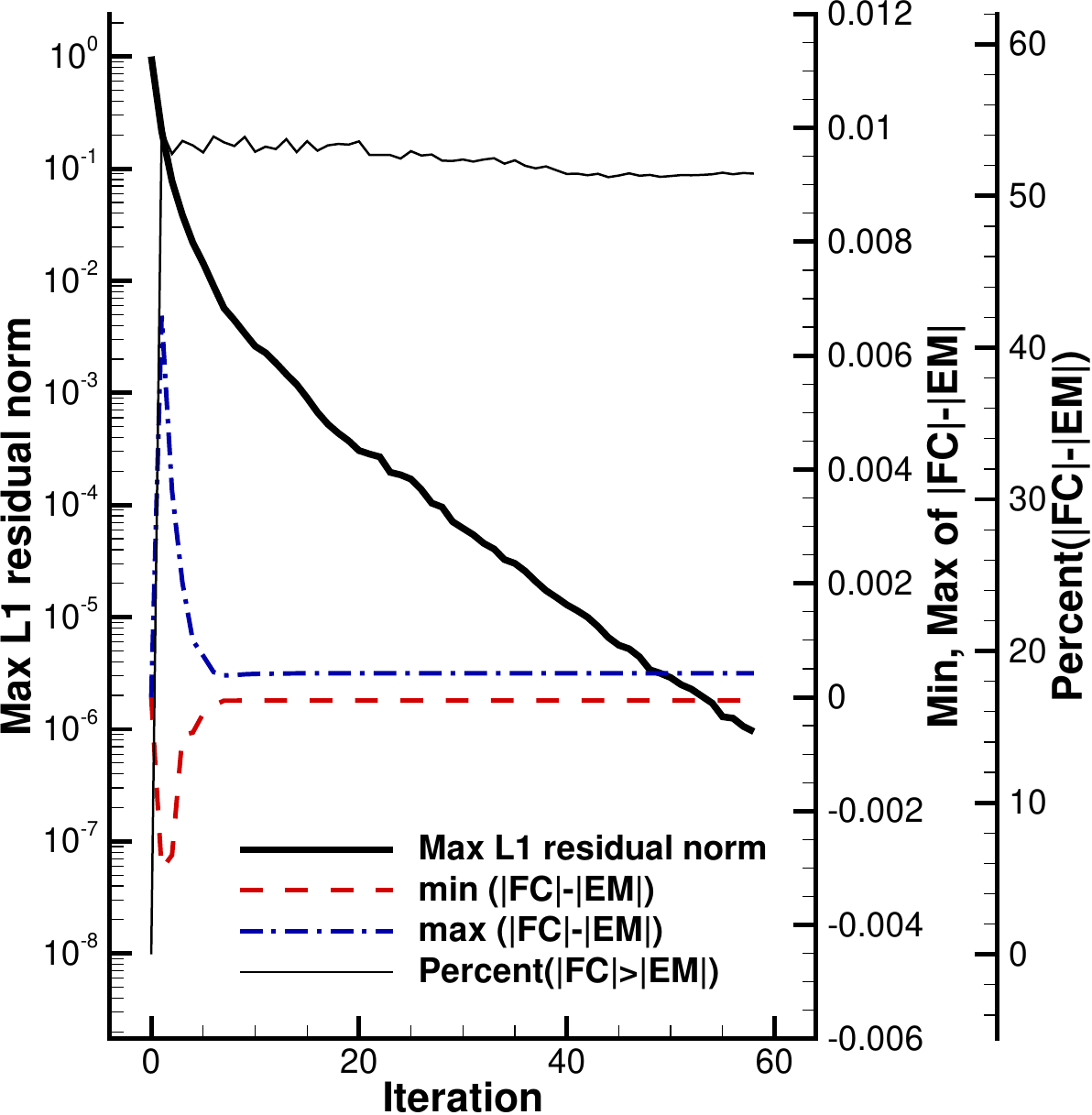}
          \caption{EM: $\alpha=4/3$, $\alpha_{Jac}=4/3$.}
          \label{fig:twod_conv_em_4o3}
      \end{subfigure}
     \hfill
          \begin{subfigure}[t]{0.32\textwidth}
    \centering
        \includegraphics[width=0.99\textwidth,trim=0 0 0 0,clip]{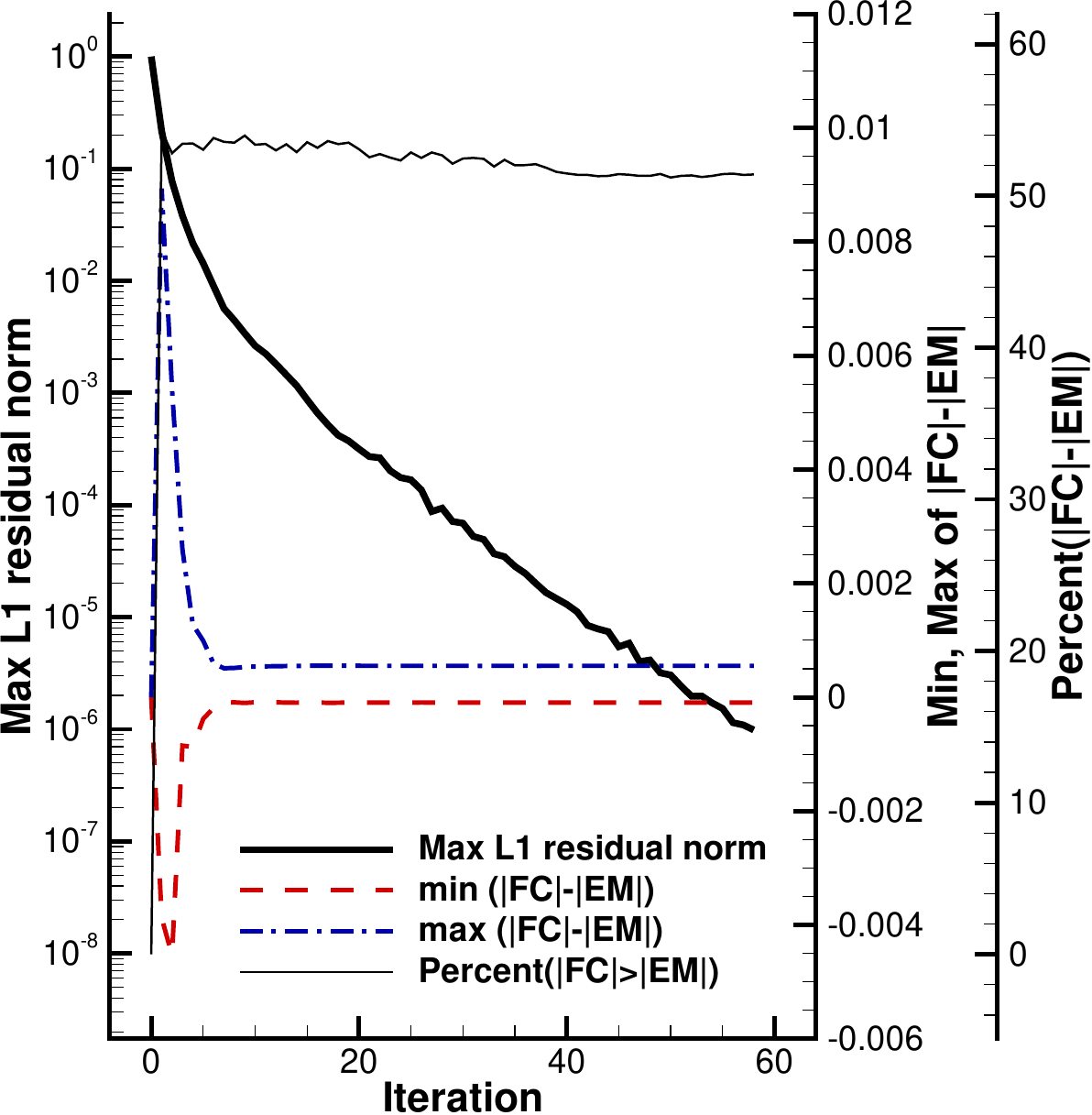}
          \caption{EM: $\alpha=1$, $\alpha_{Jac}=1$.}
          \label{fig:twod_conv_em_1p0}
      \end{subfigure} 
          \begin{subfigure}[t]{0.32\textwidth}
    \centering
        \includegraphics[width=0.99\textwidth,trim=0 0 0 0,clip]{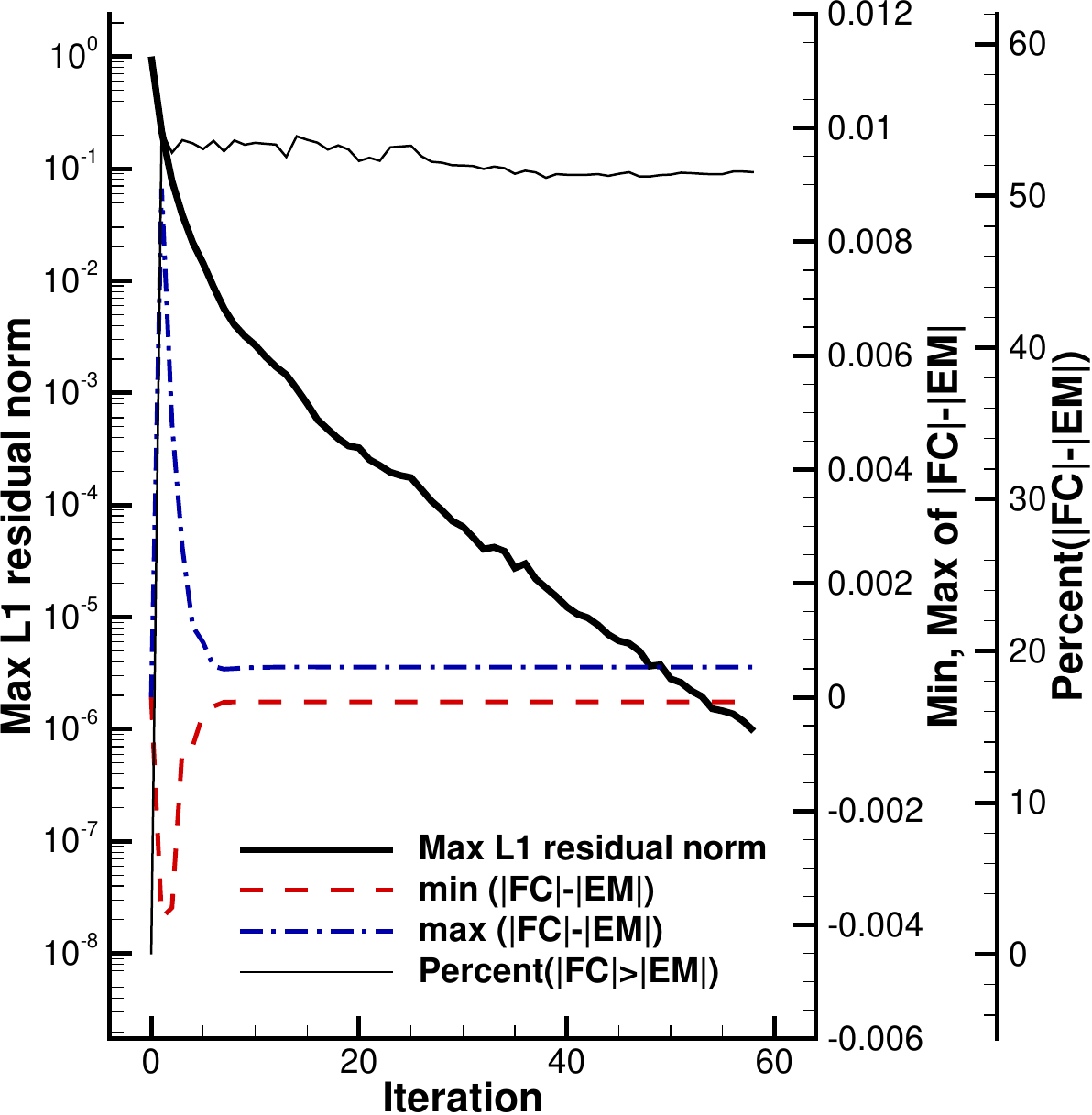}
          \caption{UR: $\alpha=4/3$, $\alpha_{Jac}=4/3$.}
          \label{fig:twod_conv_unrecon}
      \end{subfigure} 
            \caption{
          \label{fig:twod_conv}
Iterative convergence results in the flat plate case. 
}
\end{figure}

Iterative convergence results are shown in Figure \ref{fig:twod_conv}. First, Figure \ref{fig:twod_conv_fc_4o3} shows the convergence of the case of the FC scheme with $\alpha=4/3$ and $\alpha_{Jac}=4/3$. The implicit solver fails to converge; it repeats the process of divergence, CFL reduction, convergence, and then divergence again. Note that the FC scheme with $\alpha=4/3$ and $\alpha_{Jac}=4/3$ is the only case where the CFL number is changed during the iteration (CFL is kept $10^{15}$ in all other cases). Also plotted are the minimum and maximum of the differences of the magnitude of the solution jumps computed by the FC and EM schemes over all faces, and the percentage of the faces where the FC jump magnitude is greater than the EM jump. Slightly more than 50\% of the faces have the FC jump magnitude larger than the EM jump magnitude, but the FC jump magnitude is significantly larger than the EM jump magnitude (blue dashed-dotted line) as soon as the solver starts to diverge. At the same time, the EM jump magnitude is not as large when it is larger than the FC jump magnitude (red dashed line). 

Then, we tested two techniques as possible ways to stabilize the implicit solver for the FC scheme, based on the one-dimensional study. One is to increase the damping coefficient $\alpha_{Jac}$ from $4/3$ to $5$. As shown in Figure \ref{fig:twod_conv_fc_4o3_largeJac}, $\alpha_{Jac}=  5$ stabilizes the solver and successfully achieves convergence in 235 iterations. The FC jump magnitude is still very large compared with the EM jump magnitude, but it is bounded. Another is the quadratic LSQ gradient. As can be seen in Figure \ref{fig:twod_conv_fc_4o3_qlsq}, the quadratic LSQ gradient also stabilizes the solver and achieves convergence at fewer iterations (60 iterations). Observe that the FC jump magnitude is significantly reduced by the use of quadratic LSQ gradients. 

%
%

On the other hand, the solver converges rapidly in 58 iterations with the EM scheme, as shown in Figure \ref{fig:twod_conv_em_4o3}. Although the EM jump magnitude is larger than the FC jump magnitude over nearly 50\% of the faces, the maximum difference is quite small when it is larger than the FC jump magnitude (red dash line). In contrast, the FC jump magnitude is greater when it is larger than the EM jump magnitude (blue dash-dot line); this causes the solver to diverge as we have seen in Figure \ref{fig:twod_conv_fc_4o3}. 

 Also, we tested the EM scheme with $\alpha=1$, which corresponds to the classical scheme, i.e., the face-tangent scheme \cite{Haselbacher_PhD,thomas_diskin_nishikawa:CandF2011} or equivalently the Mathur-Murthy scheme \cite{MathurMurthy:NHT1997scheme}. The results are very similar to the case of $\alpha=4/3$: compare Figures \ref{fig:twod_conv_em_4o3} and \ref{fig:twod_conv_em_1p0}. The solver converged at 58 iterations as in the case of $\alpha=4/3$. This indicates that the EM scheme alone stabilizes the implicit solver (not for a particular choice of $\alpha$), effectively reducing the solution jump magnitude; and this appears to be the major reason for the classical scheme's success. 

\begin{figure}[htbp!]
\begin{center}
\begin{minipage}[b]{0.45\textwidth}
\begin{center}
\includegraphics[width=0.8\textwidth,trim=0 0 0 0,clip]{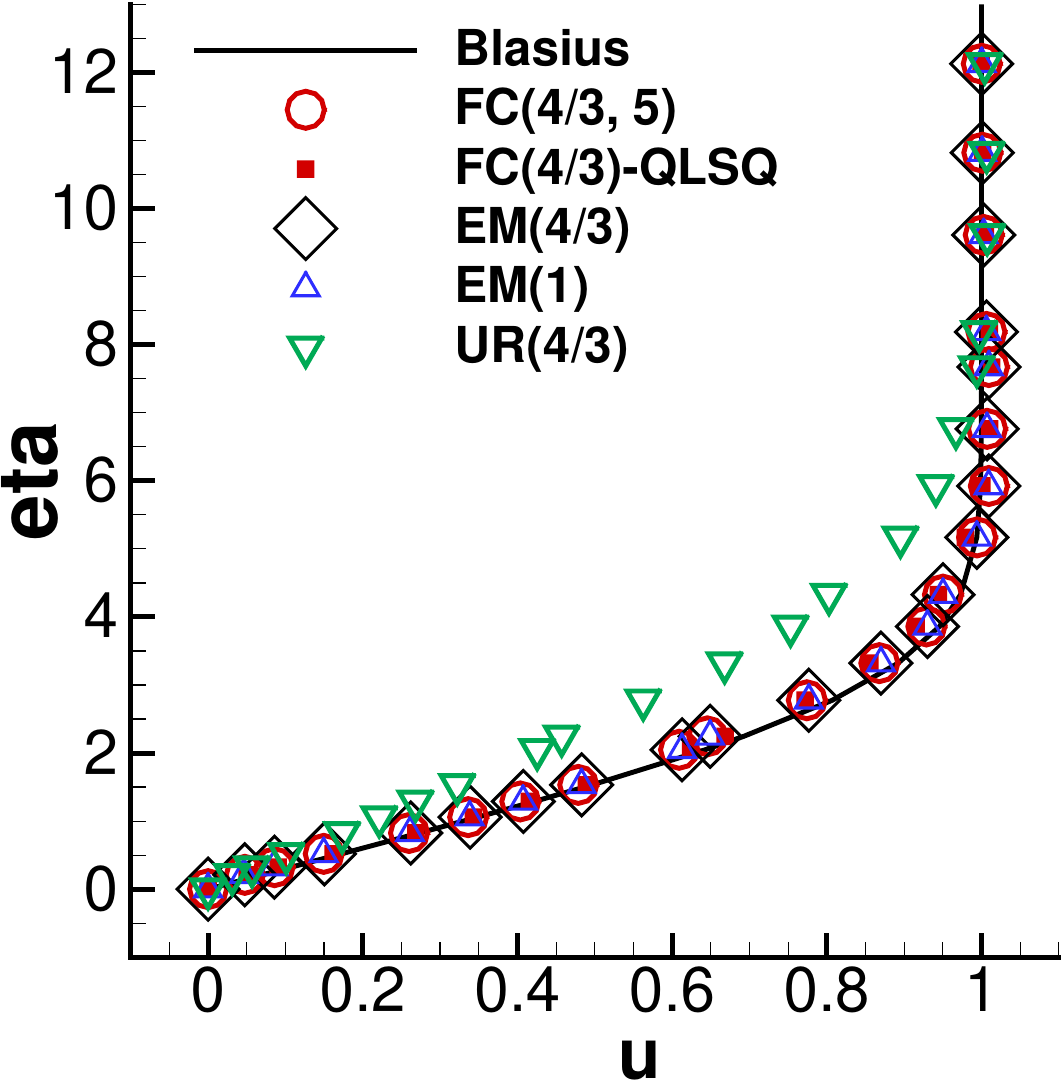}
\caption{Velocity profile comparison at $x=0.9$.}
\label{fig:twod_fp_u_profile}
\end{center}
\end{minipage}
\end{center}
\end{figure}

Finally, we tested the unreconstructed scheme, where the solution jump in the alpha-damping scheme is computed as the difference between the cell-center values, $u_R - u_L = u_k - u_j$. The solver converged successfully at 58 iterations and the convergence history is very similar to the other two EM cases, as shown in Figure \ref{fig:twod_conv_unrecon}. However, as mentioned earlier, and also demonstrated in Ref.~\cite{Nishikawa_RobustFluxes:jcp2020}, the unreconstructed solution jump makes the viscous scheme inconsistent (not even first-order accurate). This accuracy deterioration is manifested in the velocity profile sampled at $x=0.9$. See Figure \ref{fig:twod_fp_u_profile}, where numerical solutions for the converged cases are compared with the Blasius solution, and FC(4/3, 5) corresponds to the FC scheme with $\alpha=4/3$ and $\alpha_{Jac}=5$. The data with green lower triangular symbols, labeled UR(4/3),  corresponds to the unreconstructed scheme. Clearly, the solution is not accurate while others match the Blasius solution reasonably well. We may have the unreconstructed option as the last resort but should be aware of the accuracy deterioration.

\subsection{Laminar Flow over a Circular Cylinder at $Re_\infty = 20$}
\label{results_cylinedr} 

\begin{figure}[htbp!]
    \centering
          \begin{subfigure}[t]{0.48\textwidth}
    \centering
        \includegraphics[width=0.99\textwidth,trim=0 0 0 0,clip]{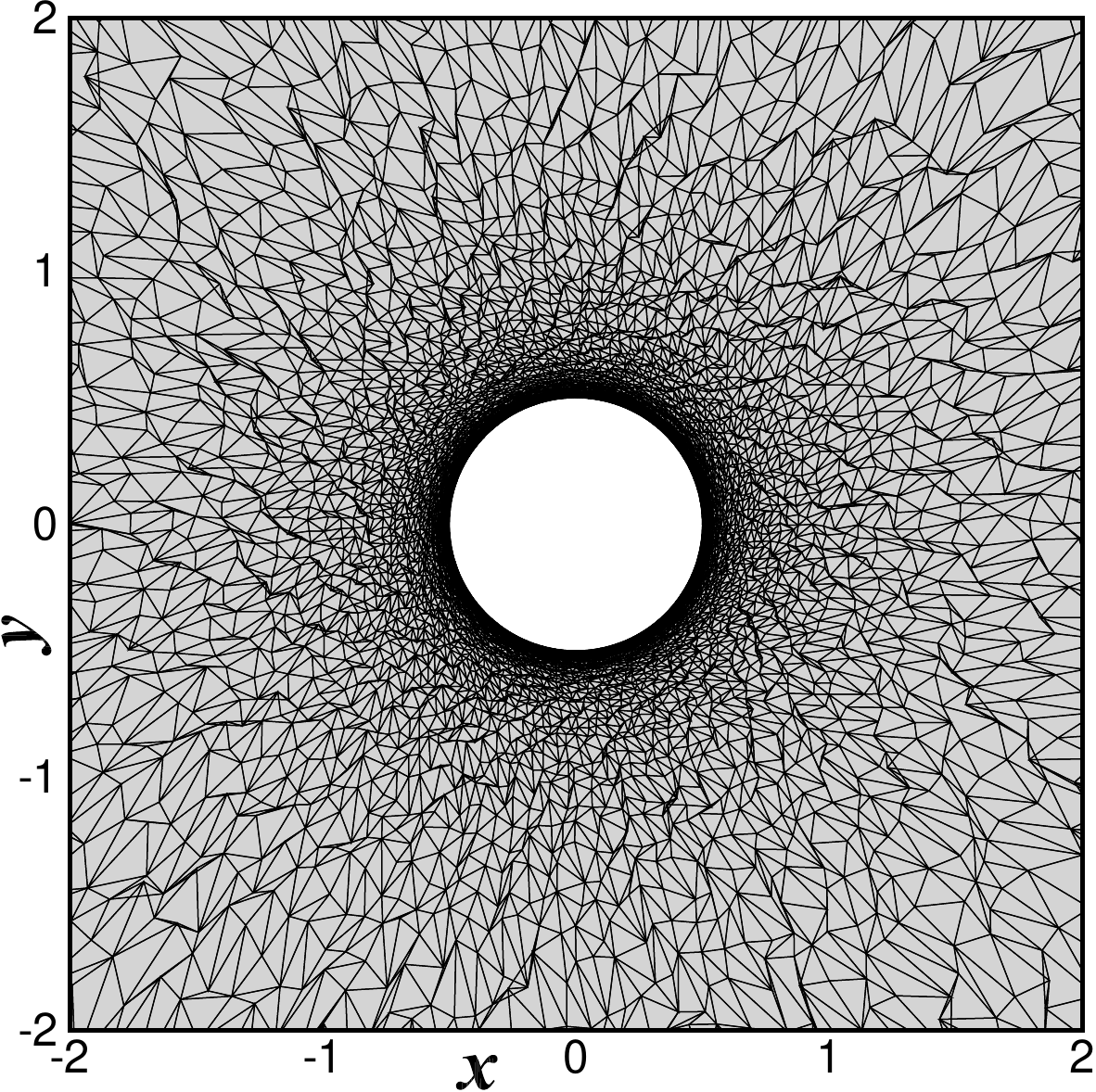}
          \caption{Irregular triangular grid.}
          \label{fig:twod_cylinder_grid}
      \end{subfigure}
      \hfill
          \begin{subfigure}[t]{0.48\textwidth}
    \centering
        \includegraphics[width=0.99\textwidth,trim=0 0 0 0,clip]{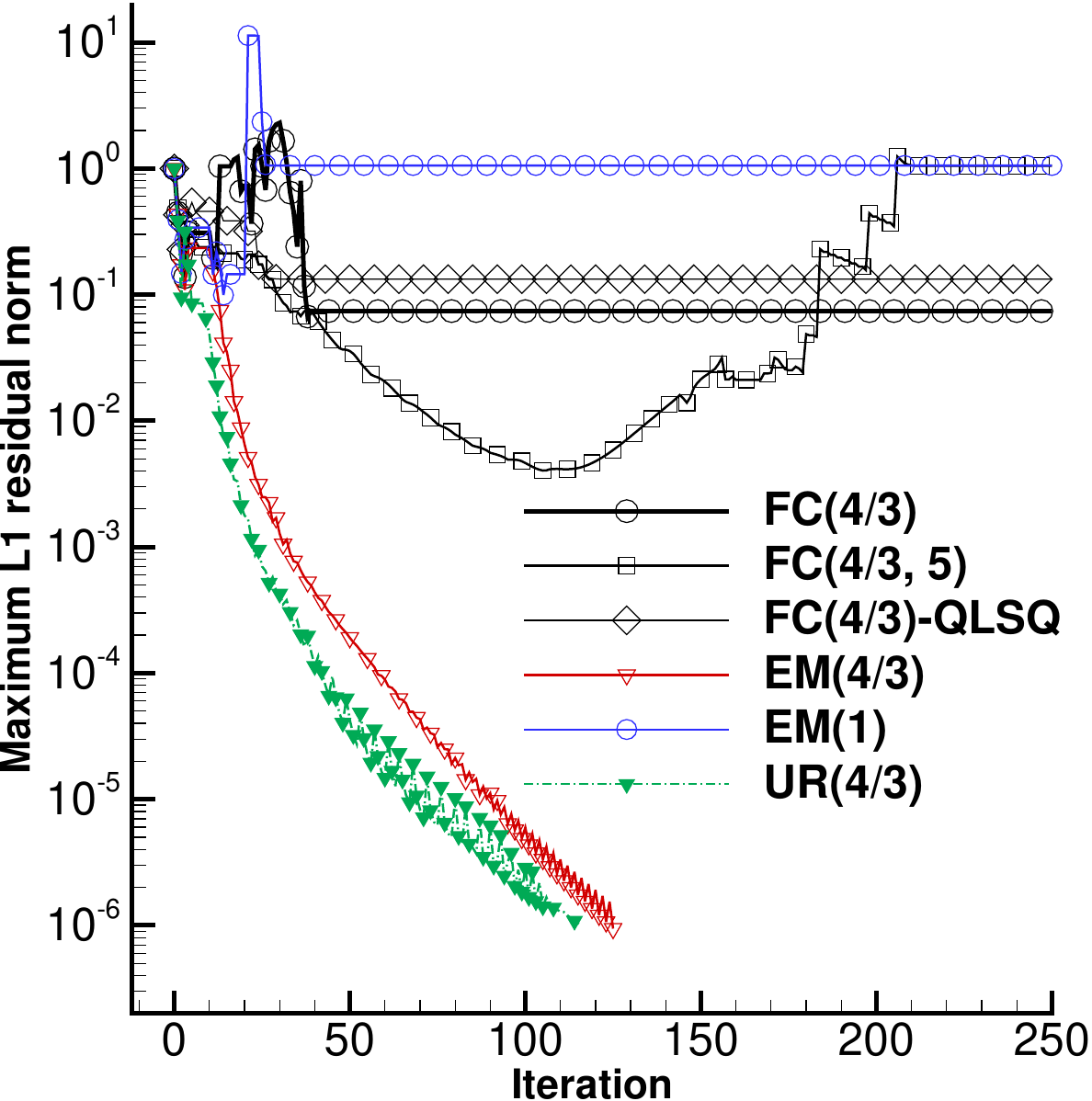}
          \caption{Iterative convergence.}
          \label{fig:twod_cylinder_convergence}
      \end{subfigure}
            \caption{
          \label{fig:twod_cylinder}
Grid and iterative convergence in the cylinder case.
}
\end{figure}

To show that the implicit solver is not always stabilized as in the previous case, we consider a steady viscous flow over a circular cylinder of unit diameter at $M_\infty = 0.15$ and $Re_\infty = 20$. The domain is defined by the cylinder centered at the origin and an outer circular boundary at the distance 50. The no-slip condition is applied at the cylinder, and a free stream condition is applied at the outer boundary. The grid is an irregular triangular grid with 103,686 nodes, 21,098 triangles, as shown in Figure \ref{fig:twod_cylinder_grid}. This grid was generated from a structured quadrilateral grid with random diagonal splitting followed by nodal perturbation. It is a highly distorted grid with some inverted (negative-volume) cells, and some CFD solvers may easily fail to converge. Note, however, that robust implicit finite-volume solvers can converge for grids with negative-volume cells \cite{nishikawa_aiaa2017-4295}. As in the previous case, we solve the residual equations by the implicit defect-correction solver with the same settings for the linear relaxation, initial values, and the CFL number.  The solver is taken to be converged when the maximum $L_1$ norm of the residuals is reduced by six orders of magnitude.

%
%

\begin{figure}[htbp!]
    \centering
          \begin{subfigure}[t]{0.48\textwidth}
    \centering
        \includegraphics[width=0.99\textwidth,trim=0 0 0 0,clip]{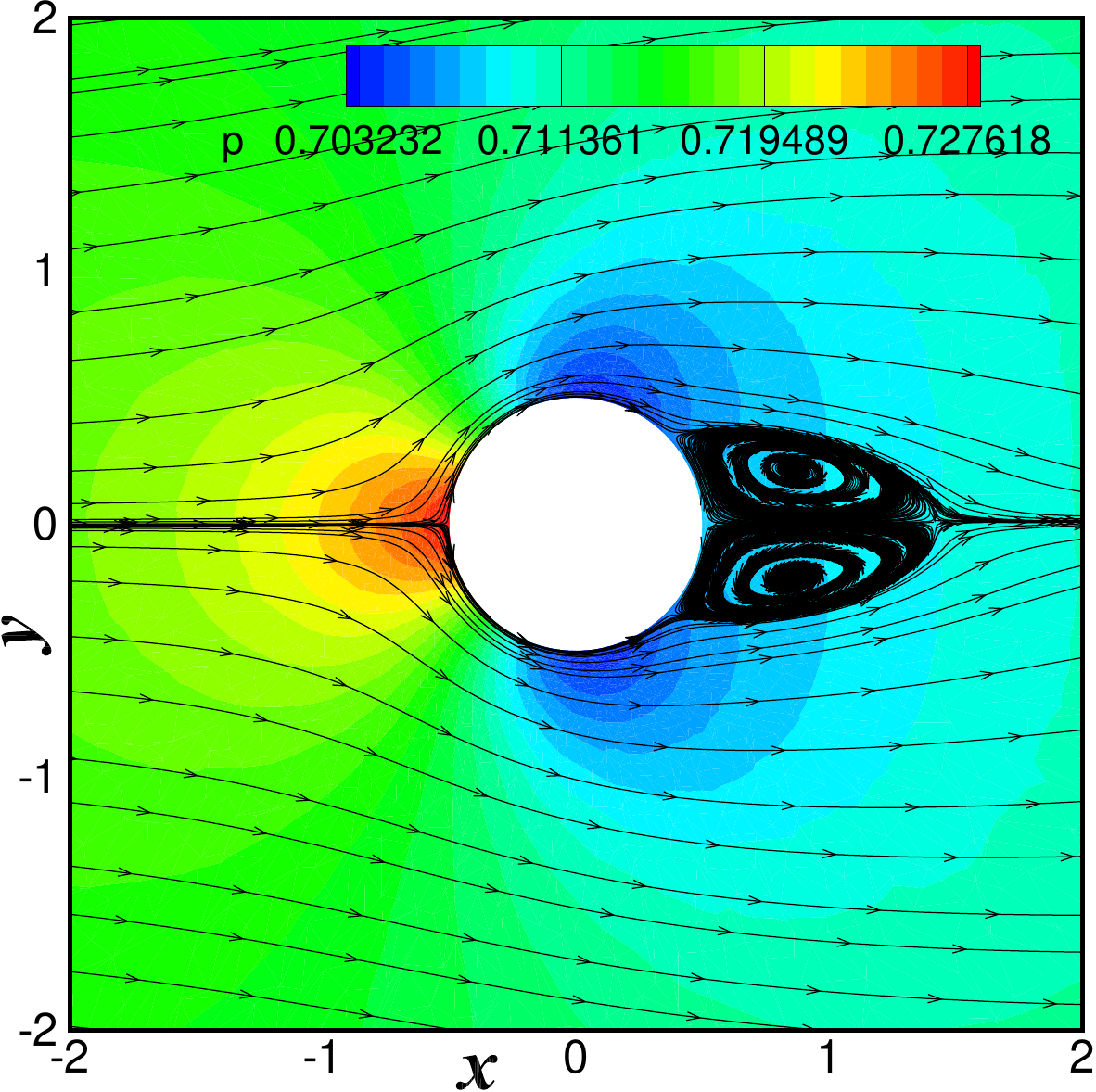}
          \caption{EM(4/3).}
          \label{fig:twod_cylinder_press_em4o3}
      \end{subfigure}
      \hfill
          \begin{subfigure}[t]{0.48\textwidth}
    \centering
        \includegraphics[width=0.99\textwidth,trim=0 0 0 0,clip]{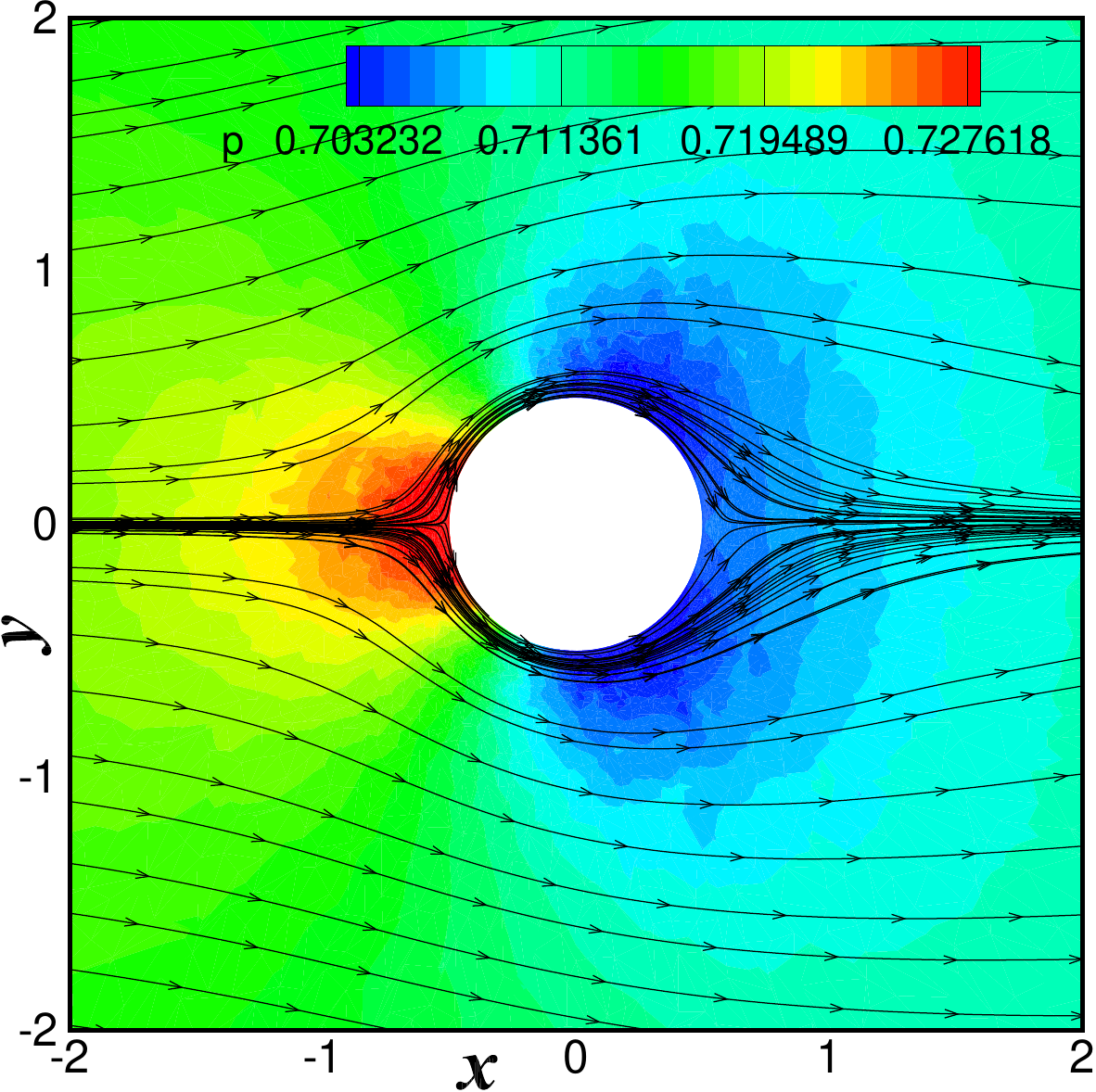}
          \caption{UR(4/3).}
          \label{fig:twod_cylinder_press_ur4o3}
      \end{subfigure}
            \caption{
          \label{fig:twod_cylinder_press}
Solutions obtained in the cylinder case.
} 
\end{figure}

Figure \ref{fig:twod_cylinder_convergence} shows the iterative convergence results. As can be seen, the implicit solver diverged for all the FC schemes, and even for the EM scheme with $\alpha =1$ (i.e., the classical scheme). The solver converged only for the remaining two schemes: the EM and unreconstructed schemes with $\alpha = 4/3$. However, the unreconstructed scheme gives, in this test case, a very inaccurate solution. Figure \ref{fig:twod_cylinder_press} compares the pressure contours and streamlines for the two converged cases. Figure \ref{fig:twod_cylinder_press_em4o3} shows the result for the EM scheme. This solution is reasonably accurate in that the length of the pair of vortices behind the cylinder is about 0.9, which is in good agreement with the value reported in the literature (see Figure 16 of Ref.~\cite{coutanceau_bouard_JFM1997}). Furthermore, the drag coefficient is found to be 1.99, which is in good agreement with the value reported in the literature \cite{Fornberg1980} (i.e., in the range 1.998-2.053). On the other hand, the numerical solution obtained with the unreconstructed scheme is very inaccurate as can be clearly seen in Figure \ref{fig:twod_cylinder_press_ur4o3}. The pressure contours are very irregular, and the streamlines show no presence of the wake vortices. Also, the drag coefficient is 2.63, which is significantly off the expected range 1.998-2.053 \cite{Fornberg1980}.

It is interesting that the EM scheme converged with $\alpha =4/3$ but not with $\alpha=1$. Note that the same value of $\alpha$ is used for the residual and the Jacobian, in both cases. As expected, increasing $\alpha$ in the Jacobian, we found that the solver converged for the residual with $\alpha=1$, producing a similar numerical solution. 

Again, the unreconstructed option may be useful for obtaining a numerical solution in cases, where any other scheme does not lead to solver convergence, but accuracy deterioration should be expected and simulation results must be carefully assessed.

%
%

\section{Conclusions}
\label{conclusions} 

\indent

In this paper, we investigated the stability of an implicit defect-correction solver widely used in practical CFD codes and demonstrated that the edge-midpoint variant of the alpha-damping viscous scheme can be more robust because of the reduced magnitude of the solution jump in the damping term. A Fourier analysis shows that the implicit solver can be made stable by increasing the damping coefficient in the Jacobian or by using a high-order solution reconstruction. Effectively, both contribute to making the magnitude of the damping term in the residual smaller than that in the Jacobian, thereby pushing down the spectral radius below 1. 

In two dimensions, we have shown that the magnitude of the solution jump computed with the edge-midpoint reconstruction scheme has a smaller upper bound than that computed with the face-centroid reconstruction scheme. As expected, two-dimensional results based on an irregular mixed-element grid show that the implicit solver is stable with the edge-midpoint reconstruction scheme. It was also shown that the implicit solver with the face-centroid reconstruction scheme is unstable but can be made stable by increasing the damping coefficient in the Jacobian and also by upgrading the accuracy of gradients from a linear least-squares method to a quadratic least-squares method. A further study using a highly distorted triangular grid has shown that the edge-midpoint reconstruction scheme with the damping coefficient $\alpha=4/3$ converged and produced an accurate numerical solution while any other schemes or techniques failed to achieve convergence. 

We have also tested an unreconstructed scheme, where the damping term is computed with the difference of the numerical solutions at two adjacent cell centers. The implicit solver converged for both test problems, but accuracy deterioration is observed because it destroys second-order accuracy of the viscous scheme, leading to inconsistency (not even first-order accurate). This option may be kept as the last resort, but one should be aware of the accuracy deterioration and carefully assess numerical results. 

It is emphasized that the edge-midpoint reconstruction scheme does not guarantee the stability of implicit solvers. The main point of the present study is that the implicit solver tends to be more stable if the damping term is made smaller relative to that in the Jacobian, and this is effectively accomplished by the edge-midpoint reconstruction scheme, thus making it a strong candidate for a default scheme. 

Future work includes extensions of the edge-midpoint reconstruction scheme to the inviscid scheme, in the framework of cell-centered finite-volume methods, for possible further improvements in robustness for irregular unstructured grids. As in the case of the viscous scheme, the solution jump in the dissipation term can be evaluated at any position, which does not affect the order of accuracy as long as it vanishes for linear solutions.

\section*{Acknowledgments} 

The author gratefully acknowledges support by 
the High Speed Flight Project, through the Hypersonic Airbreathing Propulsion Branch of NASA's Langley Research Center, funded under Contract No. 80LARC23DA003. The author is grateful to Boris Diskin for valuable discussions.


\bibliography{./bibtex_nishikawa_database}

\begin{thebibliography}{10}

\bibitem{LIU201788}
Y.~Liu and W.~Zhang.
\newblock Accuracy preserving limiter for the high-order finite volume method
  on unstructured grids.
\newblock {\em Comput. Fluids}, 149:88--99, 2017.

\bibitem{tong2026hybrid}
Y.~Tong and P.~Tsoutsanis.
\newblock A hybrid finite-volume reconstruction framework for efficient
  high-order shock-capturing on unstructured meshes.
\newblock {\em Computers \& Fluids}, 308, 2026.

\bibitem{NakashimaNishikawaLeeCerizza_aiaa_aviation2025}
Y.~Nakashima, H.~Nishikawa, J.~Lee, and D.~Cerizza.
\newblock A reduced-memory multicolor gauss-seidel relaxation scheme for
  implicit unstructured-polyhedral-grid {CFD} solver on {GPU}.
\newblock In {\em AIAA Aviation 2025 Forum}, {AIAA} Paper 2025-3870, Las Vegas,
  NV, 2025.

\bibitem{DiskinLiuNishikawa_GradFree_aiaa_aviation2026}
B.~Diskin, Y.~Liu, and H.~Nishikawa.
\newblock Efficient gradient-free linear reconstruction scheme for edge-based
  solvers.
\newblock In {\em AIAA Aviation 2026 Forum}, {AIAA} Paper 2026-4024, San Diego,
  CA, 2026.

\bibitem{Kleb_etal_aiaa2019-2948}
W.~{L}. Kleb, M.~{A}. Park, W.~{A}. Wood, K.~{L}. Bibb, K.~{B}. Thompson, and
  R.~{J}. Gomez.
\newblock Sketch-to-solution: An exploration of viscous {CFD} with automatic
  grids.
\newblock In {\em 24th {AIAA} Computational Fluid Dynamics Conference}, {AIAA}
  Paper 2019-2948, Dallas, TX, 2019.

\bibitem{ThompsonNishikawaPadway_aiaa_scitech2023}
K.~{B}. Thompson, H.~Nishikawa, and E.~Padway.
\newblock Economical third-order methods for accurate surface heating
  predictions on simplex element meshes.
\newblock In {\em AIAA SciTech 2023 Forum}, {AIAA} Paper 2023-2629, National
  Harbor, MD \& Online, 2023.

\bibitem{MoriscoNishikawa_aiaa_scitech2025-0302}
C.~{T}. Morisco and H.~Nishikawa.
\newblock Applications of implicit edge-based gradients to third-order
  edge-based scheme for adaptive tetrahedral grids.
\newblock In {\em AIAA SciTech 2025 Forum}, {AIAA} Paper 2025-0302, Orlando,
  FL, 2025.

\bibitem{Nastac2024ClosedLoop}
Gabriel~C. Nastac, Zachary Ernst, Alexandra~M. Hickey, Aaron~C. Walden,
  Kevin~E. Jacobson, William~T. Jones, Eric~J. Nielsen, Boris Diskin, Li~Wang,
  Ashley~M. Korzun, Patrick~J. Moran, Hayden~V. Dean, Bradford~E. Robertson,
  and Dimitri Mavris.
\newblock Closed-loop simulations of human-scale mars lander descent
  trajectories on frontier.
\newblock In {\em AIAA Aviation Forum and ASCEND}, Las Vegas, NV, July 2024.
\newblock AIAA 2024-3535.

\bibitem{Liu2025SlidingMesh}
Yi~Liu, Boris Diskin, Eric~J. Nielsen, Li~Wang, and Gabriel~C. Nastac.
\newblock Assessment of sliding-mesh method for predicting multi-rotor uam
  aircraft performance on graphics processing units.
\newblock In {\em AIAA Aviation Forum and ASCEND}, Las Vegas, NV, July 2025.
\newblock AIAA 2025-3511.

\bibitem{thomas_diskin_nishikawa:CandF2011}
J.~L. Thomas, B.~Diskin, and H.~Nishikawa.
\newblock A critical study of agglomerated multigrid methods for diffusion on
  highly-stretched grids.
\newblock {\em Comput. Fluids}, 41(1):82--93, February 2011.

\bibitem{nishikawa:AIAA2010}
H.~Nishikawa.
\newblock Beyond interface gradient: A general principle for constructing
  diffusion schemes.
\newblock In {\em 40th {AIAA} Fluid Dynamics Conference and Exhibit}, {AIAA}
  Paper 2010-5093, Chicago, IL, 2010.

\bibitem{Haselbacher_PhD}
A.~C. Haselbacher.
\newblock {\em A Grid-Transparent Numerical Method for Compressible Viscous
  Flow on Mixed Unstructured Meshes}.
\newblock PhD thesis, Loughborough University, 1999.

\bibitem{MathurMurthy:NHT1997scheme}
S.~{R}. Mathur and J.~{Y}. Murthy.
\newblock A pressure-based method for unstructured meshes.
\newblock {\em Numerical Heat Transfer, {P}art {B}: Fundamentals: An
  International Journal of Computation and Methodology}, 31(2):195--216, 1997.

\bibitem{nakashima_watanabe_nishikawa:Japan2014}
Y.~Nakashima, N.~Watanabe, and H.~Nishikawa.
\newblock Development of an effective implicit solver for general-purpose
  unstructured {CFD} software.
\newblock In {\em The 28th Computational Fluid Dynamics Symposium}, C08-1,
  Tokyo, Japan, 2014.

\bibitem{WhiteBaurlePasseSpiegelNishikawa:JANNAF}
J.~{A}. White, R.~Baurle, B.~{J}. Passe, S.~{C}. Spiegel, and H.~Nishikawa.
\newblock Geometrically flexible and efficient flow analysis of high speed
  vehicles via domain decomposition, part 1, unstructured-grid solver for high
  speed flows.
\newblock In {\em {JANNAF} 48th Combustion 36th Airbreathing Propulsion, 36th
  Exhaust Plume and Signatures, 30th Propulsion Systems Hazards, Joint
  Subcommittee Meeting, Programmatic and Industrial Base Meeting}, Newport
  News, VA, 2017.

\bibitem{BATISTIC2026115056}
I.~Batisti\'{c}, P.~Castrillo, and P.~Cardiff.
\newblock High-order cell-centred finite-volume solid mechanics using a
  jacobian-free newton-krylov method.
\newblock {\em Journal of Computational Physics}, 563:115056, 2026.

\bibitem{BellostaAbergoNishikawa_aiaa2025-0072}
T.~Bellosta, L.~Abergo, and H.~Nishikawa.
\newblock Enhancing accuracy in mixed-element grids and convergence on skewed
  grids for the two-dimensional edge-based compressible {N}avier-{S}tokes
  solver.
\newblock In {\em AIAA SciTech 2025 Forum}, {AIAA} Paper 2025-0072, Orlando,
  FL, 2025.

\bibitem{nakashima_private}
Yoshitaka Nakashima.
\newblock Private Communication, 2015.

\bibitem{nishikawa_nakashima_watanabe:jcp2017}
H.~Nishikawa, Y.~Nakashima, and N.~Watanabe.
\newblock Effects of high-frequency damping on iterative convergence of
  implicit viscous solver.
\newblock {\em J. Comput. Phys.}, 348:66--81, 2017.

\bibitem{nishikawa_stencil:JCP2019}
H.~Nishikawa.
\newblock Efficient gradient stencils for robust implicit finite-volume solver
  convergence on distorted grids.
\newblock {\em J. Comput. Phys.}, 386:486--501, 2019.

\bibitem{Nishikawa_RobustFluxes:jcp2020}
H.~Nishikawa.
\newblock Robust numerical fluxes for unrealizable states.
\newblock {\em J. Comput. Phys.}, 408:109244, 2020.

\bibitem{jalali_etal:CF2014}
A.~Jalali, M.~Sharbatdar, and C.~Ollivier-Gooch.
\newblock Accuracy analysis of unstructured finite volume discretization
  schemes for diffusive fluxes.
\newblock {\em Comput. Fluids}, 101:220--232, September 2014.

\bibitem{NishikawaDiskin_Note_ViscosityAveraging:2022}
H.~Nishikawa and B.~Diskin.
\newblock Arithmetic averages of viscosity coefficient are sufficient for
  second-order finite-volume viscous discretization on unstructured grids.
\newblock arXiv:2203.08334v1 [math.NA], 2022.

\bibitem{liu_nishikawa_aiaa2016-3969}
Y.~Liu and H.~Nishikawa.
\newblock Third-order inviscid and second-order hyperbolic {N}avier-{S}tokes
  solvers for three-dimensional inviscid and viscous flows.
\newblock In {\em 46th {AIAA} Fluid Dynamics Conference}, {AIAA} Paper
  2016-3969, Washington, D.C., 2016.

\bibitem{nishikawa_aiaa2017-4295}
H.~Nishikawa.
\newblock Uses of zero and negative volume elements for node-centered
  edge-based discretization.
\newblock In {\em 23rd {AIAA} Computational Fluid Dynamics Conference}, {AIAA}
  Paper 2017-4295, Denver, Colorado, 2017.

\bibitem{coutanceau_bouard_JFM1997}
M.~Coutanceau and R.~Bouard.
\newblock Experimental determination of the main features of the viscous flow
  in the wake of a circular cylinder in uniform translation. {P}art 1. {S}teady
  flow.
\newblock {\em J. Fluid Mech.}, 79:231--256, 1977.

\bibitem{Fornberg1980}
B.~Fornberg.
\newblock A numerical study of steady viscous flow past a circular cylinder.
\newblock {\em Journal of Fluid Mechanics}, 98(4):819--855, 1980.

\end{thebibliography}
\bibliographystyle{unsrt}

\end{document}